\documentclass[12pt]{amsart}

\usepackage{scrpage2}

\usepackage{Commands}
\usepackage{jabbrv}
\def\lra{\longrightarrow}

\allowdisplaybreaks

\numberwithin{equation}{section}
\title[ASIP for Random Systems]{An almost sure invariance principle for several classes of random dynamical systems}
\author{Jason Atnip}
\address{Department of Mathematics, University of North Texas, Denton, TX 76203-1430, USA}
\email{\href{jason.atnip@unt.edu}{jason.atnip@unt.edu} \hspace*{0.42cm} \it Web: \rm \href{atnipmath.com}{atnipmath.com}}

\begin{document}

%\tableofcontents
%\newpage
\begin{abstract}
	In this paper we deal with a large class of dynamical systems having a version of the spectral gap property. Our primary class of systems comes from random dynamics, but we also deal with the deterministic case.	We show that if a random dynamical system has a fiberwise spectral gap property as well as an exponential decay of correlations in the base, then, developing on Gou\"{e}zel's approach, the system satisfies the almost sure invariance principle. The result is then applied to uniformly expanding random systems like those studied by Denker and Gordin and Mayer, Skorulski, and Urba\'nski.  
\end{abstract}
\maketitle
\section{Introduction}

The almost sure invariance principle (ASIP) is a powerful statistical property which assures that the trajectories of a process can be approximated in an almost sure manner with the trajectories of a Brownian motion with a negligible error term relative to the length of the trajectories. In particular, the ASIP implies many limit theorems including the law of the iterated logarithm and various versions of the central limit theorem. For more consequences of the ASIP, see \cite{philipp_almost_1975} and the references therein.

The ASIP was first shown for scalar--valued independent and identically distributed random variables by Strassen in \cite{strassen_invariance_1964,strassen_almost_1967} and then for $\RR^d$--valued observables by Melbourne and Nicol in \cite{melbourne_vector-valued_2009}. In \cite{gouezel_almost_2010}, Gou\"ezel uses spectral properties to show a vector valued ASIP for a wide class of dynamical systems which satisfy the so called ``spectral gap" property. Concerning random dynamical systems, the ASIP has been shown for random expanding dynamical systems by Aimino, Nicol, and Vaienti in \cite{aimino_annealed_2015},considering only stationary measures. In the recent paper \cite{dragicevic_almost_2016} of Dragi\v cevi\'c, Froyland, Gonz\'alez--Tokman, and Vaienti for random Lasota--Yorke maps they consider non--stationary fiberwise measures, as we do in this paper, but they prove the ASIP for centered observables. Here we consider general H\"older observables, which are not necessarily centered.

In this paper we build upon Gou\"{e}zel's approach, to present a real valued ASIP for a large class of random dynamical systems with non--stationary fiberwise random measures for which only the central limit theorem and law of the iterated logarithm were previously known. We show that if a random dynamical system has transfer operators which satisfy the spectral gap property as well as a base dynamical system which exhibits an exponential decay of correlations, then the random system satisfies an ASIP. In particular it is the difficulty of dealing with non--centered observables which requires this extra condition on the dynamical system in the base.

In what follows we will show that the uniformly expanding random systems of \cite{mayer_distance_2011}, and in particular the DG*--systems of \cite{mayer_distance_2011} based upon the work of Denker--Gordin \cite{denker_gibbs_1999}, are examples of such well--behaved systems for which our theory applies.% the random transcendental systems of \cite{mayer_random_2014}, and shifts of finite type with weakly positive transfer operators as in \cite{mayer_countable_2015} are examples of such well behaved systems for which our theory applies.

\section{Preliminaries}~\\

\subsection{Almost Sure Invariance Principle}~\\

We will consider a real valued stationary stochastic process $(A_n)_{n=0}^\infty$ which is bounded in $L^p$ for some $p>2$. 
\begin{definition}
	Suppose $0<\lm\leq 1/2$ and $\sg^2>0$. We say that the sequence $(A_n)_{n=0}^\infty$ satisfies an almost sure invariance principle with error exponent $\lm$ and limiting covariance $\sg^2$ if there exists a probability space $\Om$ and stochastic process $(A_n')_{n=0}^\infty$ and $(B_n)_{n=0}^\infty$ on $\Om$ such that the following hold:
	\begin{enumerate}
		\item The processes $(A_n)_{n=0}^\infty$ and $(A_n')_{n=0}^\infty$ have the same distribution. 
		\item The random variables $B_n$ are independent and distributed as $\cN(0,\sg^2)$. 
		\item Almost surely in $\Om$ we have that 
		\begin{align*}
			\absval{\sum_{j=0}^{n-1}A_j'-\sum_{j=0}^{n-1}B_j}\leq o(n^\lm)
		\end{align*}
		when $n\to\infty$.
	\end{enumerate} 
	As a Brownian motion on the integers corresponds with a sum of i.i.d. Gaussian random variables, this definition can be restated as almost sure approximation by a Brownian motion. 
\end{definition}
~\\
\subsection{Random Dynamical Systems}~\\

Suppose $(X,\sF, m)$ is a complete (Borel) probability space with metric $d_X$ and that $\ta:X\lra X$ is an invertible map, often referred to as \textit{the base map}. We assume that $\ta$ preserves the measure $m$, i.e.
\begin{align*}
	m\circ\ta^{-1}=m,
\end{align*}
and that $\ta$ is ergodic with respect to $m$. For each $x\in X$ we associate the metric space $(\jl_x,\varrho_x)$ with each $x\in X$, and let 
\begin{align*}
\jl:=\union_{x\in X}\set{x}\times\jl_x%\sub X\times Y.
\end{align*}
For ease of exposition we will identify $\jl_x$ with $\set{x}\times\jl_x$.
Further suppose that for each $x\in X$ there is a continuous map $T_x:\jl_x\lra\jl_{\ta(x)}$, and define the associated skew product map $T:\jl\lra\jl$ by 
\begin{align*}
T(x,z)=(\ta(x), T_x(z)). 
\end{align*}
For each $n\geq 0$ we denote 
\begin{align*}
	T^n_x:=T_{\ta^{n-1}(x)}\circ\cdots\circ T_x:\jl_x\lra\jl_{\ta^n(x)}.
\end{align*} 
Similarly, we have 
\begin{align*}
	T^n(x,z)=(\ta^n(x),T^n_x(z)). 
\end{align*}
Given a continuous function $g:\jl\lra\RR$, for each $x\inX$ we let 
\begin{align*}
	g_x:=g\rvert_{\set{x}\times\jl_x}:\jl_x\lra\RR.
\end{align*} 

\begin{comment}
--------------------------------------------

A canonical example of such a system is to consider the finite collection $I=\set{f_1,\dots,f_m}$ together with some probability vector $p=(p_1,\dots,p_m)$. We let the base space be the full shift $\Om=I^\ZZ$, the Cartesian product of $I$ with itself $\ZZ$--many times, together with the shift map $\ta$ and product measure $m$. In this case we may think of each $\om\in\Om$ as a sequence of the $f_i$ and we then investigate the almost sure properties of the associated dynamical systems.  

----------------------------------------------

content...
\end{comment}
For a more thorough treatment of random dynamics see, for example, \cite{arnold_random_1998} or \cite{kifer_thermodynamic_2008}.
\\
\subsection{Random Measures}~\\

Suppose $T:\cJ\lra\cJ$ is a random dynamical system over the base $(X,\sF,m,\ta)$ the as defined above, and let $\sB:=\sB_\cJ$ denote the Borel $\sg$--algebra on $\cJ$ such that the following hold
\begin{enumerate}
	\item The map $T$ and projection function $\pi_X:\cJ\lra X$, given by $\pi_X(x,y)=x$, are measurable,
	\item for every $A\in\sB$, $\pi_X(A)\in\sF$,
	\item $\sB_x:=\sB\rvert_{\cJ_x}$ is a $\sg$--algebra on $\cJ_x$.
\end{enumerate}

A measure $\mu$ on $(\jl,\sB)$ is said to be \textit{random probability measure} relative to $m$ if it has marginal $m$, i.e. if
\begin{align*}
	\mu\circ\pi_X^{-1}=m.
\end{align*}
If $(\mu_x)_{x\in X}$ are disintegrations of $\mu$ with respect to the partition $\left(\cJ_x\right)_{x\in X}$ of $\cJ$, then these satisfy the following properties:
\begin{enumerate}
	\item For every $B\in\sB_x$, the map $X\ni x\longmapsto\mu_x(B)\in[0,1]$ is measurable, 
	\item For $m$-a.e. $x\in X$, the map $\sB_x\ni B\longmapsto\mu_x(B)\in[0,1]$ is a Borel probability measure.
\end{enumerate}
We let $\cP(\jl_x)$ denote the space of probability measures on $(\jl_x,\sB_x)$ for each $x\in X$, and by $\cP_m(\jl)$ we denote the space of all random measures on $\jl$ with marginal $m$. By definition we then have that for $\mu\in\cP_m(\jl)$ and suitable $g:\jl\lra\RR$
\begin{align*}
	\mu(g)=\int_{\jl}g d\mu=\int_X\int_{\jl_x}g_xd\mu_xdm(x).
\end{align*}
For such a function $g$, we denote the expected value of $g$ with respect to the measure $\mu$ by 
\begin{align*}
	\cE_\mu(g)=\mu(g).
\end{align*}
For a more detailed discussion of random measures see Crauel's book \cite{crauel_random_2002}.\\

\subsection{Continuous and H\"older Potentials}\label{sec:potentials}~\\

For each $x\in X$ we let $\cC(\cJ_x)$ be the set of all continuous and bounded functions $g_x:\cJ_x\lra\RR$. Taken together with the sup norm, $\norm{\spot}_{x,\infty}$, over the fiber $\cJ_x$, $\cC(\cJ_x)$ becomes a (fiberwise) Banach space. We can then consider the global space $\sC(\jl)$ of functions $g:\jl\lra\RR$ such that for $m$-a.e. $x\inX$ the functions $g_x:=g\rvert_{\jl_x}\in\cC(\cJ_x).$  Let $\sC^0(\jl)$ and $\sC^1(\jl)$ be subspaces of $\sC(\jl)$ such that for all $g\in\sC^0(\jl)$ the function $x\longmapsto\norm{g_x}_{x,\infty}$ is $\sF$--measurable and $g\in\sC^1(\jl)$ implies that 
\begin{align*}
	\norm{g}_1:=\int_X\norm{g_x}_{x,\infty} dm(x)<\infty.
\end{align*}
Let $\sC_*^\infty(\cJ)$ denote the space of all $\sB_\cJ$--measurable functions $g\in\sC(\jl)$ such that 
\begin{align*}
	\sup_{x\inX}\norm{g_x}_{x,\infty}<\infty.
\end{align*}
Clearly $\sC^\infty_*(\cJ)$ becomes a Banach space when coupled together with the norm $\absval{\makebox[1ex]{\textbf{$\cdot$}}}_\infty$ given by
\begin{align*}
	\absval{g}_\infty=\sup_{x\inX}\norm{g_x}_{x,\infty}, \quad g\in\sC^\infty_*(\cJ).%\int_X\norm{g_x}_{x,\infty} \,dm(x)
\end{align*}
Fix $\al\in(0,1]$. The $\al$--variation of a function $g_x\in\cC(\cJ_x)$ is given by
\begin{align*}
	v_{x,\al}(g_x)&=\sup\set{\frac{\absval{g_x(y)-g_x(y')}}{\varrho_x^\al(y,y')}:y,y'\in \cJ_x, \; 0\neq \varrho_x(y,y')\leq \eta}.
	%\\&=\inf\set{H_x:\absval{u_x(z_1)-u_x(z_2)}\leq H_x \varrho_x^\al(z_1,z_2)}.
\end{align*}
We let $\cH_\al(\cJ_x)$ denote the collection of all functions $g_x\in\cC(\cJ_x)$ such that $v_{x,\al}(g_x)<\infty$. Taken together with the norm given by
\begin{align*}
	\norm{\spot}_{x,\al}=\norm{\spot}_{x,\infty}+v_{x,\al}(\spot),
\end{align*}  
$\cH_\al(\cJ_x)$ becomes a Banach algebra, that is we have
\begin{align*}
	\norm{g_xh_x}_{x,\al}\leq \norm{g_x}_{x,\al}\cdot\norm{h_x}_{x,\al}	
\end{align*}
for $g_x,h_x\in\cH_\al(\cJ_x)$. We say that a function $g\in\sC^1(\cJ)$ is (global) $\al$--H\"older continuous over $\cJ$ if there is a $m$--measurable function
\begin{align*}
	H:X\lra[1,\infty), \qquad x\longmapsto H_x,
\end{align*}   
such that $v_{x,\al}(g_x)\leq H_x$ for $m$-a.e. $x\in X$ and $\log H\in L^1(m)$, i.e.
\begin{align*}
	\int_X \log H_x \,dm(x)<\infty.
\end{align*}
Let $\sH_\al(\jl)$ be the collection of all such functions. For each $1\leq p<\infty$ we let $\sH_\al^p(\jl)$ denote the set of functions $g\in\sH_\al(\jl)$ such that 
$\norm{g_x}_{x,\al}\in L^p(m)$, that is such that 
\begin{align*}
	\trinorm{g}_{\al,p}:=\left(\int_X\norm{g_x}_{x,\al}^p\,dm(x)\right)^{1/p}<\infty.
\end{align*}
We then have that $\sH_\al^p(\jl)$ taken with the norm $\trinorm{\spot}_{\al,p}$ is a Banach space. 
Now let $\sH_\al^*(\jl)$ be the set of all functions $g\in\sC^\infty_*(\cJ)$ such that \begin{align*}
	\sup_{x\inX}v_{x,\al}(g_x)<\infty.
\end{align*} 
For $g\in\sH_\al^*(\cJ)$, we set 
\begin{align*}
	V_\al(g)=\sup_{x\inX}v_{x,\al}(g_x), %\int_X v_{x,\al}(g_x)\, dm(x),
\end{align*}
and then define the norm $\absval{\makebox[1ex]{\textbf{$\cdot$}}}_\al$ on $\sH_\al^*(\cJ)$ by 
\begin{align*}
	\absval{g}_\al=\absval{g}_\infty+V_\al(g).
\end{align*}
Then $\sH_\al^*(\jl)$ is a Banach space when considered with the norm $\absval{\makebox[1ex]{\textbf{$\cdot$}}}_\al$. 

We also consider the set $\HH_\bt(X)$, for $\bt\in(0,1]$, of bounded continuous functions $G:X\lra\RR$ such that the $\bt$--variation of $G$, denoted by $\vta_\bt(G)$, is finite, i.e. 
\begin{align*}
	\vta_\bt(G):=\sup\set{\frac{\absval{G(x)-G(x')}}{d_X^\bt(x,x')}:x,x'\in X, \; x\neq x'}<\infty.
\end{align*}
We say that such functions are $\bt$--H\"older continuous on $X$. $\HH_\bt$ then becomes a Banach space, and in fact a Banach algebra, when coupled together with the norm
\begin{align*}
	\norm{\spot}_{\HH}:=\norm{\spot}_X+\vartheta_\bt(\spot),
\end{align*}
where $\norm{\spot}_X$ denotes the supremum norm on $X$.\\

%Following in Gou\"ezel's notation, for the purposes of simplification, if $a,b\in\RR$ we let $e^{iab}$ denote $e^{i\langle a,b\rangle}$, where $\langle a,b\rangle$ is the usual inner product on $\RR$. 

\subsection{Transfer Operators}~\\
Given a continuous functions $g\in\sC^1(\cJ)$, for each $n\in\NN$ and $x\in X$ we then define the Birkhoff sum $S_{x,n}:\cC(\cJ_x)\lra\RR$ by 
\begin{align*}
S_{x,n}g_x:=\sum_{j=0}^{n-1}g_{\ta^j(x)}\circ T_x^j.
\end{align*}
If there is no confusion about the fiber $\cJ_x$, we will simply write $S_n$.
Now given a function $\phi\in\sH_\al(\cJ)$ we define the (Perron--Frobenius) transfer operator $\tr_{\phi,x}:\cC(\cJ_x)\lra\cC(\cJ_{\ta(x)})$ by 
\begin{align*}
\tr_x(u_x)(w):=\tr_{\phi, x}(u_x)(w):=\sum_{z\in T_x^{-1}(w)}u_x(z)e^{\phi_x(z)}, \quad w\in\cJ_{\ta(x)}.
\end{align*}
For $n\in\NN$ the iterates $\tr^n_x:\cC(\cJ_x)\lra\cC(\cJ_{\ta^n(x)})$ is given by 
\begin{align*}
	\tr_x^n(u_x):=\tr_{\ta^{n-1}(x)}\circ\dots\circ\tr_x(u_x).
\end{align*}
Inductively one can show that 
\begin{align*}
	\tr_x^n(u_x)(w)=\sum_{z\in T^{-n}_x(w)}u_x(z)e^{S_n\phi_x(z)}, \quad w\in\cJ_{\ta^n(x)}
\end{align*}
for each $n\in\NN$. 
Denote by $\tr_x^*$ the dual operator $\tr_x^*:\cC^*(\cJ_{\ta(x)})\lra\cC^*(\cJ_x)$, where $\cC^*(\cJ_x)$ is the dual space of $\cC(\cJ_{x})$ equipped with the weak$^*$ topology. Now suppose there is a random probability measure $\nu$ on $\cJ$ such that 
\begin{align*}
\tr_x^*\nu_{\ta(x)}=\lm_x\nu_x \quad \text{ for } m-\text{a.e. }x\inX, 
\end{align*}
where
\begin{align*}
\lm_x:=\tr_x^*(\nu_{\ta(x)})(\ind)=\nu_{\ta(x)}(\tr_x\ind)=\int_{\cJ_{\ta(x)}}\tr_x(\ind_x)\,d\nu_{\ta(x)}.
\end{align*}
We are then able to define the normalized operator $\tr_{0,x}:\cC(\cJ_x)\lra\cC(\cJ_{\ta(x)})$ by 
\begin{align*}
	\tr_{0,x}(u_x):=\lm_x^{-1}\tr_x(u_x).
\end{align*}
Clearly we have that the iterates $\tr^n_{0,x}:\cC(\cJ_x)\lra\cC(\cJ_{\ta^n(x)})$ of the normalized are given by 
\begin{align*}
	\tr_{0,x}^n(u_x)=(\lm_{x}^n)^{-1}\tr_x^n(u_x).
\end{align*}
where 
\begin{align*}
		(\lm_{x}^n)^{-1}=(\lm_{x})^{-1}\cdots(\lm_{\ta^{n-1}(x)})^{-1}=\int_{\cJ_{\ta^n(x)}}\tr_x^n(\ind_x)\,d\nu_{\ta^n(x)}
\end{align*} 
For each $r\in\RR$ we define the perturbed operator given by
\begin{align*}
\tr_{r,x}(u_x):=\tr_{0,x}(e^{irg_x}\cdot u_x)
\end{align*} 
In the sequel, the perturbed operator will be our main technical tool.
The following lemma characterizes the iterates of the perturbed operator.
\begin{lemma}\label{lem: Lrn...Lr0=L0(stuff)}
	For $r_0,\dots, r_{n-1}\in\RR$ we have 
	\begin{align*}
	\tr_{r_{n-1},\ta^{n-1}(x)}\circ\cdots\circ\tr_{r_0,x}(u_x)&=\tr^n_{0,x}(e^{i\sum_{j=0}^{n-1}r_jg_{\ta^j(x)}\circ T^j_x}\cdot u_x).
	\end{align*}
\end{lemma}
\begin{proof}
	For $r_0,r_1\in\RR$ we have 
	\begin{align*}
	\tr_{r_1,\ta(x)}\left(\tr_{r_0,x}(u_x)\right)
	&=\tr_{0,\ta(x)}(e^{ir_1g_{\ta(x)}}\cdot\tr_{0,x}(e^{ir_0g_x}\cdot u_x))\\
	&=\tr_{0,\ta(x)}\left(\tr_{0,x}\left( e^{ir_1g_{\ta(x)}\circ T_x}\cdot e^{ir_0g_x}\cdot u_x\right)\right)\\
	&=\tr_{0,x}^2\left(e^{i(r_0g_x+r_1g_{\ta(x)}\circ T_x)}\cdot u_x\right)
	\end{align*}
	Inducting on $n\geq 1$ we suppose that for $r_{n-1},\dots,r_0\in\RR$ we have
	\begin{align*}
	\tr_{r_{n-1},\ta^{n-1}(x)}\circ\cdots\circ\tr_{r_0,x}(u_x)&=\tr^n_{0,x}(e^{i\sum_{j=0}^{n-1}r_jg_{\ta^j(x)}\circ T^j_x}\cdot u_x),
	\end{align*} 
	and let $r_n\in\RR$.  Then 
	\begin{align*}
	\tr_{r_n,\ta^n(x)}\tr_{r_{n-1},\ta^{n-1}(x)}\circ\cdots\circ\tr_{r_0,x}(u_x)&=\tr_{r_n,\ta^n(x)}\tr^n_{0,x}(e^{i\sum_{j=0}^{n-1}r_jg_{\ta^j(x)}\circ T^j_x}\cdot u_x)\\	
	&=\tr_{0,\ta^n(x)}\left(e^{ir_ng_{\ta^n(x)}}\cdot\tr^n_{0,x}(e^{i\sum_{j=0}^{n-1}r_jg_{\ta^j(x)}\circ T^j_x}\cdot u_x)\right)\\
	&=\tr^{n+1}_{0,x}(e^{i\sum_{j=0}^{n}r_jg_{\ta^j(x)}\circ T^j_x}\cdot u_x)
	\end{align*}

\end{proof}

In Section \ref{Sec: Rand G} we present our a main result, a theorem which establishes an almost sure invariance principle for quite general classes of random dynamical systems. Afterward we provide several classes of examples for which our theorem applies. Throughout our paper, $C$ will denote some positive constant, which may change from line to line.  By $\cN(\mu,\sg^2)$ we mean the normal distribution with mean $\mu$ and variance $\sg^2$.

\section{ASIP for Random Systems}\label{Sec: Rand G}
In this section we give the main result of the paper, which is an adaptation of  Theorem \ref{thm: gouezel thm} (Theorem 2.1 of \cite{gouezel_almost_2010}) for random dynamical systems. We will follow the general strategy of Gou\"ezel's proof.
Our main theorem is the following.
\begin{theorem}\label{thm: modified gouezel asip}	
	Suppose $(X,\cF,m,\ta)$ and $T:\jl\lra\jl$ form a random dynamical system as defined above and suppose that $\mu$ is a $T$--invariant measure on $\cJ$. For $\al\in(0,1]$ and let $g, \phi\in\sH_\al^*(\cJ)$. Suppose the transfer operators $\tr_{\phi,x}$, $\tr_{0,x}$, and $\tr_{r,x}$ are defined as above and suppose there is $\ep_0>0$, $\rho_x\in\cH_\al(\cJ_x)$, and $\nu_x\in\cP(\cJ_x)$ for each $x\in X$ such that the following hold.
	\begin{enumerate}
		\item \label{assum 1}There exists $C>0$ such that for each $x\inX$ 
		\begin{align*}
		\mu_x=\rho_x\nu_x \spand \norm{\rho_x}_{x,\al}\leq C.
		\end{align*} 
		Moreover, for any $r_0,\dots,r_{n-1}\in\RR$ with $\absval{r_j}<\ep_0$ we have
		\begin{align*}
		\cE_\mu\left(e^{i\sum_{j=0}^{n-1}r_jg\circ T^j}\right)=\int_X\int_{\jl_x}\tr_{r_{n-1},\ta^{-1}(x)}\circ\cdots\circ\tr_{r_0,\ta^{-n}(x)}(\rho_{\ta^{-n}(x)})\,d\nu_x\, dm(x).
		\end{align*}
		\begin{comment}
		content...
		\item \label{assum 2} For each $x\in X$ and all $u_x\in \cB_x$ we have $\cB_x\sub\cG_x\sub L^1(\nu_x)$	with 
		\begin{align*}
		\norm{u_x}_{L^1(\nu_x)}\leq\norm{u_x}_{\cG_x}\leq\norm{u_x}_{\cB_x},
		\end{align*}
		such that 
		\begin{align*}
		\int_X\norm{u_x}_{\cG_x}dm(x)<\infty \spand \int_X\norm{u_x}_{\cB_x}dm(x)<\infty .
		\end{align*}
		
		\end{comment}
		\item \label{assum 3} There exists $C\geq1$ such that for all $n\in\NN$, all $r\in\RR$ with $\absval{r}<\epsilon_0$, $m$-a.e. $x\inX$, all $f_x\in \cC(\cJ_x)$, and all $h_x\in\cH_\al(\cJ_x)$ 
		\begin{align*}
		\norm{\tr_{r,x}^nf_x}_{\ta^n(x),\infty}\leq C\norm{f_x}_{x,\infty} \spand \norm{\tr_{r,x}^nh_x}_{\ta^n(x),\al}\leq C\norm{h_x}_{x,\al}.
		\end{align*}		
		\item \label{assum 4} For each $x\in X$ there is an operator $Q_x:\cC(\cJ_x)\lra\cC(\cJ_{\ta(x)})$ defined by 
		\begin{align*}
		Q_xu_x:=\int_{\jl_x}u_x\,d\nu_x\cdot\rho_{\ta(x)},
		\end{align*} 
		and there are $C>0$ and $\kp\in(0,1)$ such that for $m$--a.e. $x\in X$ and $u_x\in\cH_\al(\cJ_x)$ 
		\begin{align*}
		\norm{\tr^n_{0,x}u_x-Q^n_xu_x}_{\ta^n(x),\infty}\leq C\kappa^n\norm{u_x}_{x,\al}.
		\end{align*}
		\item \label{assum 5} There exist constants $C>0$ and $\kp\in(0,1)$ such that for all $F\in\HH_\bt$, all $G\in L^1(m)$, and $n\in\NN$ sufficiently large we have that
		\begin{align*}
		\absval{m(G\circ \ta^{-n}\cdot F)-m(G)\cdot m(F)}\leq C\kp^n\norm{F}_\HH \norm{G}_{L^1(m)}.
		\end{align*}
		%\item We also require that $\norm{g_x}_{\cB_x}\in L^p(m)$ for some $p>2$.
		\item \label{assum 6}Suppose there exists $C>0$ and $\bt\in(0,1]$ such that for each $n\in\NN$ and $r_0,\dots,r_{n-1}\in\RR$ with $\absval{r_j}<\ep_0$ for each $1\leq j\leq n-1$, we have that the function 
		\begin{align*}
			X\ni x\longmapsto F(x)=\int_{\cJ_x}\tr_{r_{n-1},\ta^{n-1}(x)}\circ\cdots\circ\tr_{r_0,x}(\rho_x)\,d\nu_x
		\end{align*}
		is $\bt$--H\"older continuous on $X$ and furthermore, we have that 
		\begin{align*}
			\norm{F}_\HH\leq C%\norm{\rho_{\ta^{-n}(x)}}_{\bt_{\ta^{-n}}(x)}
		\end{align*} 
		independent of the choice of $n$ and $r_{0},\dots,r_{n-1}$. 
	\end{enumerate}  
	Then either there exists a real number $\sg^2>0$ such that the stochastic process $\set{g\circ T^n-\mu(g)}_{n\in\NN}$, considered with respect to the measure $\mu$, satisfies an ASIP with limiting covariance $\sg^2$ for any error exponent larger than $1/4$. Consequently, the sequence
	\begin{align*}
	\frac{S_ng-n\cdot\mu(g)}{\sqrt{n}}
	\end{align*} 
	converges in probability to $\cN(0,\sg^2)$ and 
	\begin{align*}
		\limty{n} \frac{S_ng-n\cdot\mu(g)}{\sqrt{2n\log\log n}}=1.
	\end{align*}
	Or, if $\sg^2=0$, then we have that
	\begin{align*}
		\sup_{n\in\NN}\norm{S_ng-\mu(g)}_{L^2(\mu)}<\infty 
	\end{align*}
	or equivalently, $g$ is of the form $g=k-k\circ T+\mu(g)$ where $k\in L^2(\mu)$. In addition, 
	almost surely we have
	\begin{align*}
		\limty{n}\frac{S_ng}{n^{1/4}}=0.
	\end{align*}	 
\end{theorem}
\begin{remark}
	Note that, using induction, we have that the operator
	\begin{align*}
	Q^n_x:\cC(\cJ_x)\lra\cC(\cJ_{\ta^n(x)})
	\end{align*}
	is defined by 
	\begin{align*}
	Q^n_xu_x:=\int_{\jl_x}u_x\,d\nu_x\cdot\rho_{\ta^n(x)}.
	\end{align*}
	\begin{comment}
	We also note that the requirement that 
	\begin{align*}
	\int_X\norm{u_x}_{\cG_x}dm(x)<\infty \spand \int_X\norm{u_x}_{\cB_x}dm(x)<\infty 
	\end{align*}
	could be changed depending on the context. In the case of the random shift maps in section six, we require instead that the essential supremum over $X$ is bounded rather than the integral.
	
	content...
	\end{comment}
\end{remark}

\begin{remark}\label{rem: encoding}
	We also note that in \cite{gouezel_almost_2010}, Gou\"ezel refers to hypothesis \eqref{assum 1} above as ``\textit{the characteristic function of the process $g\circ T^n$ being encoded by the family of operators $\set{\tr_r}$}". 
\end{remark}	
\begin{remark}
	Assumption \eqref{assum 5} essentially says that the dynamical system in the base exhibits an exponential decay of correlations. While this may seem like a strong assumption, it should, in some sense, be expected. Consider a potential $\phi$ which is constant on fibers, that is, $\phi_x=c_x$ for some constant $c_x\in\RR$, for each $x\inX$. Such a potential only captures the dynamics of the base map $\ta$, and, furthermore, the ASIP simply does not hold in general assuming only ergodicity. 
	
	Also note that since $\ta$ is invertible and $m$--invariant, we also have that
	\begin{align}\label{pos exp dec corr}
	\absval{m(G\cdot F\circ\ta^n)-m(G)\cdot m(F)}\leq C\kp^n\norm{G}_\HH \norm{F}_{L^1(m)},
	\end{align}
	for $n\in\NN$ and $F,G\in \HH_\bt\cap L^1(m)$.
	
	As we will see in the sequel, there are many random systems which satisfy this assumption. In particular, any random system with an expanding base map $\ta$ and Gibbs measure $m$ will satisfy an exponential decay of correlations for H\"older continuous observables, see, for example, \cite{przytycki_conformal_2010}.
\end{remark}
In order to prove Theorem \ref{thm: modified gouezel asip} we will adapt the method of the proof of the following theorem of Gou\"{e}zel:
\begin{theorem}[Theorem 2.1, \cite{gouezel_almost_2010}]\label{thm: gouezel thm}
	Let $(A_\ell)$ be a stochastic process whose characteristic function is encoded by a family of operators $\tr_r:\cB\lra\cB$ (see Remark \ref{rem: encoding}) and which is bounded in $L^p$ for some $p>2$. Further assume
	\begin{enumerate}
		\item[(I0)] There exists $u_0\in\cB$ and $\xi_0\in\cB^*$, the dual of $\cB$, such that for any $r_0,\dots,r_{n-1}\in\RR^d$ with $\absval{r_j}\leq\epsilon_0$,
		\begin{align*}
		\cE(e^{i\sum_{j=0}^{n-1}r_jA_j)})=\langle\xi_0,\tr_{r_{n-1}}\circ\cdots\circ\tr_{r_0}u_0\rangle.
		\end{align*}
		\item[(I1)] -- One can write $\tr_0=Q+S$, where $Q$ is a one--dimensional projection and $S$ is an  operator on $\cB$ with $SQ=QS=0$ and $\norm{\tr_0^n-Q}_{\cB\lra\cB}\leq C\kappa^n$ for some $\kappa<1$. 
		\item[(I2)] -- There exists $C>0$ such that $\norm{\tr_r^n}_{\cB}\leq C$ for all $n\in\NN$ and all small enough $r\in\RR^d$. 
	\end{enumerate}  
	Then there exists $a\in\RR^d$ and a matrix $\Sg^2$ such that $(\sum_{j=0}^{n-1}A_j-na)/\sqrt{n}$ converges to $\cN(0,\Sg^2)$. Moreover the process $(A_j-a)_{j\in\NN}$ satisfies an ASIP with limiting covariance $\Sg^2$ for any error exponent larger that $p/(4p-4)$.
	
\end{theorem} 

\begin{remark}	
%	We remark that the theorem of Gou\"ezel is the special case our Theorem \ref{thm: modified gouezel asip} when the base space $X$ is a singleton and the Banach space $\cG=\cB$. 
	
	We note that we assume in the first half of \eqref{assum 4} in the hypotheses of our Theorem \ref{thm: modified gouezel asip} what Gou\"{e}zel proves in his first step. In practice, we will always have such an operator, so our assumption is justified.
\end{remark}
However, in order to prove Theorem \ref{thm: modified gouezel asip}, we must invoke Gou\"ezel's Theorem 1.3 and Lemma 2.7 of \cite{gouezel_almost_2010}. First we state the main assumption of \cite{gouezel_almost_2010}, which ensures that a given stochastic process is sufficiently close to an independent process for our purposes. We will refer to this assumption as condition (H).

\begin{enumerate}
	\item[\mylabel{cond H}{(H)}] There exists $\ep_0>0$ and constants $C,c>0$ such that for any $n,m,k>0$, $b_1<b_2<\dots<b_{n+m+1}$, and $r_1,\dots, r_{n+m}\in\RR^d$ with $\absval{r_j}\leq \ep_0$ we have
	\begin{align*}
	\bigg\lvert \cE & \left(e^{i\sum_{j=1}^n t_j(\sum_{\ell=b_j}^{b_{j+1}-1}A_\ell)+i\sum_{j=n+1}^{n+m}t_j(\sum_{\ell=b_j+k}^{b_{j+1}+k-1}A_\ell)}\right) \\
	&\left.-\cE\left(e^{i\sum_{j=1}^nt_j(\sum_{\ell=b_j}^{b_{j+1}-1}A_\ell)}\right)\cdot
	\cE\left(e^{i\sum_{j=n+1}^{n+m}t_j(\sum_{\ell=b_j+k}^{b_{j+1}+k-1}A_\ell)}\right)
	\right\rvert\\
	&\leq C\left(1+\max\absval{b_{j+1}-b_j}\right)^{C(n+m)}e^{-ck}.
	\end{align*}
\end{enumerate}

\begin{theorem}[Theorem 1.3, \cite{gouezel_almost_2010}]\label{thm: gouezel thm 1.3}
	Let $(A_0,A_1, \dots)$ be a centered $\RR^d$--valued process, bounded in $L^p$ for some $p>2$, satisfying condition \ref{cond H}. Assume, moreover, that $\sum\absval{A_\ell}<\infty$ and that there exists a matrix $\Sg^2$ such that, for any $\al>0$, 
	\begin{align*}
	\absval{\cov\left(\sum_{\ell=m}^{n+m-1}A_\ell\right)-n\Sg^2}\leq Cn^\al,
	\end{align*}
	uniformly in $n,m$. The sequence $\sum_{\ell=0}^{n-1}A_\ell/\sqrt{n}$ converges in distribution to $\cN(0,\Sg^2)$. Moreover, the process $(A_0, A_1,\dots)$ satisfies an ASIP with limiting covariance $\Sg^2$ for any error exponent $\lm>p/(4p-4)$. 
\end{theorem}
\begin{lemma}[Lemma 2.7, \cite{gouezel_almost_2010}]\label{lem: gouezel lem 2.7}
	Let $(A_\ell)$ be a process bounded in $L^p$ for some $p>2$, satisfying condition \ref{cond H} such that for any $m\in\NN$ there exists a matrix $s_m$ such that uniformly in $\ell, m$ we have 
	\begin{align*}
	\absval{\cov(A_\ell,A_{\ell+m})-s_m}\leq Ce^{-\dl\ell}.
	\end{align*}
	Then the series $\Sg^2=s_0+\sum_{m=1}^\infty(s_m+s_m^*)$ converges in norm and we have that 
	\begin{align*}
	\absval{\cov\left(\sum_{\ell=m}^{m+n-1}A_\ell\right)-n\sg^2}\leq C.
	\end{align*}
\end{lemma}

\begin{proof}[Proof (of Theorem \ref{thm: modified gouezel asip}).]
	
	First we show that condition \ref{cond H} holds. Let $b_1<\dots<b_{n+m+1}$, $r_1,\dots, r_{n+m}\in\RR$, $x\in X$, and $k\in\NN$. Letting $y_x=\ta^{-(b_{n+m+1}+k)}(x)$, $w_x=\ta^{-(b_{n+m+1}-b_{n+1})}(x)$, $z_x=\ta^{-(b_{n+m+1}-b_{n+1}+k)}(x)$, and somewhat ignoring the fiberwise subscript notation for the moment, we see  
	\begin{align}
	&\cE_\mu\left(e^{i\sum_{j=1}^n r_j\left(\sum_{\ell=b_j}^{b_{j+1}-1}g\circ f^\ell\right)+i\sum_{j=n+1}^{n+m}r_j\left(\sum_{\ell=b_j+k}^{b_{j+1}+k-1} g\circ f^\ell\right) }\right)\nonumber\\
	&=\int_X\int_{\jl_x}\tr_{r_{n+m}}^{b_{n+m+1}-b_{n+m}}\circ\cdots\circ\tr_{r_{n+1}}^{b_{n+2}-b_{n+1}}\tr_0^k\tr_{r_{n}}^{b_{n+1}-b_{n}}\circ\cdots\circ\tr_{r_1}^{b_2-b_1}\tr_0^{b_1}\rho_{y_x} \,d\nu_x\,dm(x)\nonumber\\
	&=\int_X\int_{\jl_x}\tr_{r_{n+m}}^{b_{n+m+1}-b_{n+m}}\circ\cdots\circ\tr_{r_{n+1}}^{b_{n+2}-b_{n+1}}(\tr_0^k-Q^k)\tr_{r_{n}}^{b_{n+1}-b_{n}}\circ\cdots\circ\tr_{r_1}^{b_2-b_1}\tr_0^{b_1}\rho_{y_x}\, d\nu_x\, dm(x)\nonumber\\
	&+\int_X\int_{\jl_x}\tr_{r_{n+m}}^{b_{n+m+1}-b_{n+m}}\circ\cdots\circ\tr_{r_{n+1}}^{b_{n+2}-b_{n+1}}Q^k\tr_{r_{n}}^{b_{n+1}-b_{n}}\circ\cdots\circ\tr_{r_1}^{b_2-b_1}\tr_0^{b_1}\rho_{y_x}\, d\nu_x\, dm(x)\nonumber\\
	&=\int_X\int_{\jl_x}\tr_{r_{n+m}}^{b_{n+m+1}-b_{n+m}}\circ\cdots\circ\tr_{r_{n+1}}^{b_{n+2}-b_{n+1}}(\tr_0^k-Q^k)\tr_{r_{n}}^{b_{n+1}-b_{n}}\circ\cdots\circ\tr_{r_1}^{b_2-b_1}\tr_0^{b_1}\rho_{y_x} \,d\nu_x \,dm(x)\label{ineq: long int}\\
	&\qquad+\int_X\bigg[\int_{\jl_x} \tr_{r_{n+m}}^{b_{n+m+1}-b_{n+m}}\circ\cdots\circ\tr_{r_{n+1}}^{b_{n+2}-b_{n+1}} \rho_{w_x}\, d\nu_x \nonumber\\
	&\qquad\qquad\qquad\quad\cdot\int_{\jl_{z_x}} \tr_{r_{n}}^{b_{n+1}-b_{n}}\circ\cdots\circ\tr_{r_1}^{b_2-b_1}\tr_0^{b_1}\rho_{y_x} \,d\nu_{z_x}\bigg]\, dm(x)\nonumber.
	\end{align}
	Define the functions $F$ and $G$, which depend on the constants $n,m,r_1,\dots,r_{n+m}, b_1,\dots,b_{n+m+1}$, by 
	\begin{align*}
	F(x)&=\int_{\jl_x} \tr_{r_{n+m}}^{b_{n+m+1}-b_{n+m}}\circ\cdots\circ\tr_{r_{n+1}}^{b_{n+2}-b_{n+1}} \rho_{w_x} \,d\nu_x \\
	G(x)&=\int_{\jl_{w_x}} \tr_{r_{n}}^{b_{n+1}-b_{n}}\circ\cdots\circ\tr_{r_1}^{b_2-b_1}\tr_0^{b_1}\rho_{\ta^{-b_{n+m+1}}(x)}\, d\nu_{w_x} .
	\end{align*}

	Now as $z_x=\ta^{-(b_{n+m+1}-b_{n+1}+k)}(x)=\ta^{-(b_{n+m+1}-b_{n+1})}\circ\ta^{-k}(x)=w_x\circ\ta^{-k}(x)$ and $y_x=\ta^{-(b_{n+m+1}+k)}(x)=\ta^{-b_{n+m+1}}\circ\ta^{-k}(x)$ we can rewrite the second term in the above product as 
	\begin{align*}
	\int_{\jl_{z_x}}\tr_{r_{n}}^{b_{n+1}-b_{n}}\circ\cdots\circ\tr_{r_1}^{b_2-b_1}\tr_0^{b_1}\rho_{y_x}\, d\nu_{z_x}=G\circ \ta^{-k}(x),
	\end{align*}
	and thus the global integral of the products, the last two lines of the above string of equalities beginning with \eqref{ineq: long int},  becomes
	\begin{align*}
	\int_X F(x)\cdot G\circ\ta^{-k}(x)\,dm(x).
	\end{align*}
	By virtue of assumption \eqref{assum 5} we then see that 
	\begin{align*}
	\absval{\int_X F(x)\cdot G\circ\ta^{-k}(x)dm(x)-\int_XF(x)\,dm(x)\cdot\int_XG(x)\,dm(x)}\leq C\kp^k\norm{F}_{\HH}\cdot \norm{G}_{L^1(m)}.
	\end{align*}
	By assumption \eqref{assum 6} we have that $\norm{F}_{\HH}\leq C$ and we are able to estimate $\norm{G}_{L^1(m)}$ as follows. 
	\begin{align*}
	\norm{G}_{L^1(m)}&=\int_X\absval{\int_{\jl_{z_x}}\tr_{r_{n}}^{b_{n+1}-b_{n}}\circ\cdots\circ\tr_{r_1}^{b_2-b_1}\tr_0^{b_1}\rho_{y_x} \,d\nu_{z_x}}\,dm(x)\\
	&\leq \int_X\norm{\tr_{r_{n}}^{b_{n+1}-b_{n}}\circ\cdots\circ\tr_{r_1}^{b_2-b_1}\tr_0^{b_1}\rho_{y_x}}_{L^1(\nu_{w_x})} \,dm(x)\\
	&\leq \int_X \norm{\tr_{r_{n}}^{b_{n+1}-b_{n}}\circ\cdots\circ\tr_{r_1}^{b_2-b_1}\tr_0^{b_1}\rho_{y_x}}_{w_x,\infty} \,dm(x)\\
	&\leq C^{n+1}\int_X \norm{\rho_x}_{x,\al}\,dm(x).		
	\end{align*}
	Now we wish to estimate the $\cG_x$ norm of \eqref{ineq: long int}. In light of assumptions \eqref{assum 3} and \eqref{assum 4} we see 
	\begin{align*}
	&\norm{\tr_{r_{n+m}}^{b_{n+m+1}-b_{n+m}}\circ\cdots\circ\tr_{r_{n+1}}^{b_{n+2}-b_{n+1}}(\tr_0^k-Q^k)\tr_{r_{n}}^{b_{n+1}-b_{n}}\circ\cdots\circ\tr_{r_1}^{b_2-b_1}\tr_0^{b_1}\rho_{y_x}}_{x,\infty}\\
	&\leq C^{m}\norm{(\tr_0^k-Q^k)\tr_{r_{n}}^{b_{n+1}-b_{n}}\circ\cdots\circ\tr_{r_1}^{b_2-b_1}\tr_0^{b_1}\rho_{y_x}}_{w_x,\infty}\\
	&\leq C^{m+1}\kappa^k\norm{\tr_{r_{n}}^{b_{n+1}-b_{n}}\circ\cdots\circ\tr_{r_1}^{b_2-b_1}\tr_0^{b_1}\rho_{y_x}}_{z_x,\al}\\
	&\leq C^{m+n+1}\kappa^k\norm{\rho_{y_x}}_{y_x,\al}.
	\end{align*}
	Integrating over $X$ provides the global inequality
	\begin{align*}
	&\absval{ \int_X \int_{\jl_x}\tr_{r_{n+m}}^{b_{n+m+1}-b_{n+m}}\circ\cdots\circ\tr_{r_{n+1}}^{b_{n+2}-b_{n+1}}(\tr_0^k-Q^k)\tr_{r_{n}}^{b_{n+1}-b_{n}}\circ\cdots\circ\tr_{r_1}^{b_2-b_1}\tr_0^{b_1}\rho_{y_x}\, d\nu_x\, dm(x) }\\
	&\leq C^{m+n+1}\kp^k\int_X\norm{\rho_{x}}_{x,\al}\,dm(x).
	\end{align*}
	Thus combining our estimates we then have 
	\begin{align*}
	\bigg\rvert\cE_\mu&\left(e^{i\sum_{j=1}^n r_j\left(\sum_{\ell=b_j}^{b_{j+1}-1}g\circ f^\ell\right)+i\sum_{j=n+1}^{n+m}r_j\left(\sum_{\ell=b_j+k}^{b_{j+1}+k-1} g\circ f^\ell\right) }\right)\\
	&-\cE_\mu\left(e^{i\sum_{j=1}^n r_j\left(\sum_{\ell=b_j}^{b_{j+1}-1}g\circ f^\ell\right)}\right)
	\cdot\cE_\mu\left(e^{i\sum_{j=n+1}^{n+m} r_j\left(\sum_{\ell=b_j}^{b_{j+1}+k-1}g\circ f^\ell\right)}\right)\bigg\rvert\\
	&\qquad+C^{m+n+1}\kp^k\int_X\norm{\rho_{x}}_{x,\al}\,dm(x)+C^{n+3}\kp^k\cdot \int_X\norm{\rho_{x}}_{x,\al}\,dm(x) .
	\end{align*} 
	Thus the difference in condition \ref{cond H} is bounded by $C^{m+n+1}\kp^k\int_X\norm{\rho_{x}}_{x,\al}\,dm(x)$ for some $C>1$ and $0<\kp<1$. Setting $\kp=e^{-c}$, $C=2^{C'}$ then we see that 
	\begin{align*}
	C^{m+n+1}\kp^k\int_X\norm{\rho_{x}}_{x,\al}\,dm(x)\leq 2^{C'(n+m+1)}e^{-ck}L
	\leq L\cdot 2^{C'}(1+\max\absval{b_{j+1}-b_j})^{C'(n+m)}e^{-ck},
	\end{align*}
	where $L$ is some constant such that $\int_X\norm{\rho_x}_{x,\al}\,dm(x)\leq L$. This verifies that condition \ref{cond H} is satisfied.
	
	Next, Gou\"ezel shows that there exists $a\in\RR$ and $C,\dl>0$ such that 
	\begin{align}\label{gouezel step 3}
	\absval{\cE_\mu(g\circ T^n)-a}\leq Ce^{-\dl n}.
	\end{align}
	To accomplish this, he shows that the sequence $(\cE_\mu(g\circ T^n))_{n\geq 0}$ is Cauchy and must have some limit which he denotes by $a$. However, by the assumed $T$--invariance of $\mu$, we immediately have that $a$ must be  $\int_\cJ g\,d\mu$ and that the quantity \eqref{gouezel step 3} must in fact be equal to zero.
	For the next step in the proof of Theorem \ref{thm: modified gouezel asip} we claim that for any $m\in\NN$ there is a real number $s_m$ such that
	\begin{align*}
	\absval{\cov(g\circ T^n, g\circ T^{n+m})-s_m}\leq Ce^{-\dl n}
	\end{align*}
	uniformly in $n,m$. To show this, we show that $\cov(g\circ T^n, g\circ T^{n+m})$ is a Cauchy sequence in $n$ and therefore must converge to some limit, which we will call $s_m$.
	By the $T$--invariance of $\mu$ we have 
	\begin{align*}
	\cov(g\circ T^n, g\circ T^{n+m})&=\absval{\int_{\cJ}(g\circ T^n)(g\circ T^{n+m})d\mu-\int_{\cJ}g\circ T^n\,d\mu\cdot \int_{\cJ}g\circ T^{n+m}\,d\mu}\\
	&=\absval{\int_{\cJ}(g\cdot g\circ T^{m})\circ T^n d\mu-\left(\int_{\cJ}g\,d\mu\right)^2}\\
	&=\absval{\int_{\cJ}g\cdot g\circ T^{m}\,d\mu-\left(\int_{\cJ}g\,d\mu\right)^2}\\
	&=\cov(g,g\circ T^m).
	\end{align*}
	Thus to show that $\cov(g\circ T^n, g\circ T^{n+m})$ is Cauchy it suffices to show that $\left(\cE_\mu(g\cdot g\circ T^m)\right)_{m\geq 0}$ is Cauchy.
	\begin{align*}
	&\cov(g,g\circ T^m)\\
	&=\absval{\int_\cJ g\cdot g\circ T^m \,d\mu-\int_{\cJ}g \,d\mu\cdot\int_{\cJ}g\circ T^m \,d\mu}\\
	&=\absval{\int_X\int_{\cJ_x}g_x\cdot g_{\ta^m(x)}\circ T^m_x\, d\mu_x \,dm(x)-\int_X\int_{\cJ_x}g_x\,d\mu_x\,dm(x)\cdot\int_X\int_{\cJ_x}g_{\ta^m(x)}\circ T_x^m\,d\mu_x\,dm(x)}\\
	&\leq \absval{\int_X\int_{\cJ_x}g_x\cdot g_{\ta^m(x)}\circ T^m_x\, d\mu_x \,dm(x)-\int_X\left(\int_{\cJ_x}g_x\,d\mu_x\cdot\int_{\cJ_x}g_{\ta^m(x)}\circ T_x^m\,d\mu_x\right)\,dm(x)}\\
	&\quad + \left|\int_X\left(\int_{\cJ_x}g_x\,d\mu_x\cdot\int_{\cJ_x}g_{\ta^m(x)}\circ T_x^m\,d\mu_x\right)\,dm(x)\right.\\
		&\left.\qquad\quad-\int_X\int_{\cJ_x}g_x\,d\mu_x\,dm(x)\cdot\int_X\int_{\cJ_{\ta^m(x)}}g_{\ta^m(x)}\,d\mu_{\ta^m(x)}\,dm(x)\right|		
	\end{align*}
	We denote the first summand to the right of the inequality by $\Sg_1$ and the second by $\Sg_2$. Beginning with $\Sg_2$, we set $G(x)=\int_{\cJ_x}g_x\,d\mu_x$ and apply \eqref{pos exp dec corr} to obtain the estimate 
	\begin{align*}
	\Sg_2&= \absval{\int_X G(x)\cdot G\circ \ta^m(x)\,dm(x)-\int_X G(x)\,dm(x)\cdot\int_X G\circ\ta^m(x)\,dm(x)}\\
	&\leq C\kp^m\norm{G}_{\HH}\cdot\norm{G}_{L^1(m)}.
	\end{align*}
	Now setting $\hat{g}_x:=g_x-\int_{\cJ_x}g_x\, d\mu_x$ we can estimate the first summand, $\Sg_1$, as 
	\begin{align*}
	\Sg_1&=\absval{\int_X\int_{\cJ_x}g_x\cdot g_{\ta^m(x)}\circ T^m_x \,d\mu_x\, dm(x)-\int_X\left(\int_{\cJ_x}g_x\,d\mu_x\cdot\int_{\cJ_x}g_{\ta^m(x)}\circ T_x^m\,d\mu_x\right)\,dm(x)}\\
	&\leq \int_X\absval{\int_{\cJ_x}g_x\cdot g_{\ta^m(x)}\circ T^m_x \,d\mu_x -\int_{\cJ_x}\left(g_{\ta^m(x)}\circ T^m_x\right)\cdot\left(\int_{\cJ_x} g_x\,d\mu_x\right)\,d\mu_x }\,dm(x)\\
	&=\int_X\absval{\int_{\cJ_x} \left(g_{\ta^m(x)}\circ T^m_x\right)\cdot\left(g_x-\int_{\cJ_x}g_x\, d\mu_x\right)\,d\mu_x}\,dm(x)\\
	&=\int_X\absval{\int_{\cJ_x} \left(g_{\ta^m(x)}\circ T^m_x\right)\cdot\hat{g}_x\,d\mu_x}\,dm(x)\\
	&=\int_X\absval{\int_{\cJ_{\ta^m(x)}}\tr_{0,x}^m\left(\hat{g}_x\rho_x\cdot g_{\ta^m(x)}\circ T^m_x\right)\,d\nu_{\ta^m(x)}}\,dm(x)\\
	&\leq \int_X\int_{\cJ_{\ta^m(x)}}\absval{g_{\ta^m(x)}\cdot\tr_{0,x}^m(\hat{g}_x\rho_x)}\,d\nu_{\ta^m(x)}\,dm(x)\\
	&=\int_X\norm{g_{\ta^m(x)}\cdot\tr_{0,x}^m(\hat{g}_x\rho_x)}_{L^1(\nu_{\ta^m(x)})}\,dm(x)\\
	&\leq\int_X\norm{g_{\ta^m(x)}\cdot\tr_{0,x}^m(\hat{g}_x\rho_x)}_{\ta^m(x),\infty}\, dm(x)\\\
	&\leq \int_X\norm{g_{\ta^m(x)}}_{\ta^m(x),\infty}\cdot\norm{\tr_{0,x}^m(\hat{g}_x\rho_x)}_{\ta^m(x),\infty}\, dm(x).
	\end{align*}
	Now since $\int_{\cJ_x}\hat{g}_x \,d\mu_x=\int_{\cJ_x}\hat{g}_x\rho_x\,d\nu_x=0$, assumption \eqref{assum 4} gives 
	\begin{align*}
	\norm{\tr_{0,x}^m(\hat{g}_x\rho_x)}_{\ta^m(x),\infty}
	&=\norm{\tr_{0,x}^m(\hat{g}_x\rho_x)-\int_{\cJ_x}\hat{g}_x\rho_x\,d\nu_x\cdot\rho_{\ta^m(x)} }_{\ta^m(x),\infty}\\
	&\leq C\kp^m\norm{\hat{g}_x\rho_x}_{x,\al}.
	\end{align*}
	Using this we can continue to get
	\begin{align*}
	\Sg_1&\leq C\kp^m\int_X\norm{g_{\ta^m(x)}}_{\ta^m(x),\al}\cdot\norm{\hat{g}_x\rho_x}_{x,\al} \,dm(x)\\
	&\leq C\kp^m\int_X\norm{g_{\ta^m(x)}}_{\ta^m(x),\al}\cdot\norm{\hat{g}_x}_{x,\al} \,dm(x).
	\end{align*}
	Noting that 
	\begin{align*}
	\norm{\hat{g}_x}_{x,\al}\leq \norm{g_x}_{x,\al}+\norm{\int_{\cJ_x}g_x\,d\mu_x\cdot\ind_{\cJ_x}}_{x,\al}\leq2\norm{g_x}_{x,\al},
	\end{align*}
	and since $\norm{g_x}_{x,\al}\in L^p(m)$ for some $p>2$ we apply the Cauchy-Schwarz inequality to get 
	\begin{align*}
	\Sg_1&\leq C\kp^m\left(\int_X\norm{g_{\ta^m(x)}}^2_{\ta^m(x),\al}\, dm(x)\right)^{1/2}\left(\int_X\norm{\hat{g}_x}^2_{x,\al} \,dm(x)\right)^{1/2}\\
	&\leq C\kp^m\left(\int_X\norm{g_x}^2_{x,\al}\, dm(x)\right)^{1/2}\left(\int_X\norm{g_x}^2_{x,\al} \,dm(x)\right)^{1/2}\\
	&\leq C\kp^m \int_X\norm{g_x}_{x,\al}^p\,dm(x).
	\end{align*}

	Combining these two estimates we have that
	\begin{align*}
	\cov(g,g\circ T^m)\leq C\kp^m\left(\int_X\norm{g_x}_{x,\al}^p\,dm(x)+\norm{G}_{\HH}\cdot\norm{G}_{L^1(m)}\right),
	\end{align*}
	which finishes the claim.
	To finish the proof of Theorem \ref{thm: modified gouezel asip} we invoke Lemma \ref{lem: gouezel lem 2.7} and Theorem \ref{thm: gouezel thm 1.3}.
	
\end{proof}
In what follows we present several examples of random dynamical systems for which Theorem \ref{thm: modified gouezel asip} can be applied. Specifically, we give examples of systems for which assumptions \eqref{assum 1}--\eqref{assum 4} are known, or can be easily checked, in the literature. We also provide examples of base systems for which assumption \eqref{assum 5} is well known. Assumption \eqref{assum 6} will need to be checked for most systems as it requires a connection between the random fiber system as well as the system in the base, however, we shall present examples of systems under which this condition is met. 

\section{Random Distance Expanding Maps}\label{sec: RDEM}
In this section we give an overview of uniformly expanding random systems as they are defined by Mayer, Skorulski, and Urba\'nski in \cite{mayer_distance_2011}. 
Suppose $(X,\scr{B},m,\ta)$ is a measure preserving dynamical system with an invertible and ergodic map $\ta:X\lra X$. For each $x\inX$ we associate the compact metric space $(\jl_x,\varrho_x)$, which has been normalized in size such that $\diam_{\varrho_x}(\jl_x)\leq 1$. Given a $z\in\jl_x$ and $r>0$ we denote the ball of radius $r$ centered at $z$ in $(\jl_x,\varrho_x)$ by $B_x(z,r)$. Define the space $\jl$ by 
\begin{align*}
\jl=\union_{x\inX}\set{x}\times\jl_x.
\end{align*}
A map $T:\jl\lra\jl$ is called an \textit{expanding random map} if the mappings $T_x:\jl_x\lra\jl_{\ta(x)}$ are continuous open surjections and if there exists a function $\eta:X\lra\RR$, $x\longmapsto\eta_x$,  and $\xi>0$ such that the following hold:
\begin{itemize}
	\item \textit{Uniform Openness:} $T_x(B_x(z,\eta_x))\bus B_{\ta(x)}(T_x(z,\xi))$ for every $(x,z)\in\jl$.
	\item \textit{Measurably Expanding:} There exists a measurable function $\gamma:X\lra(1,\infty)$, $x\longmapsto\gamma_x$, such that for $m$-a.e. $x\in X$
	\begin{align*}
	\varrho_{\ta(x)}(T_x(z_1), T_x(z_2)\geq\gamma_x\varrho_x(z_1,z_2)\quad \text{ whenever }\quad\varrho_x(z_1,z_2)<\eta_x, \quad z_1,z_2\in\jl_x.
	\end{align*}
	\item \textit{Measurability of the Degree:} The map $x\longmapsto\deg(T_x):=\sup_{y\in\jl_{\ta(x)}}\# T_x^{-1}(\set{y})$ is measurable. 
	\item \textit{Topological Exactness:} There exists a measurable function $x\longmapsto n_\xi(x)$ such that for almost every $x\in X$ and every $z\in\jl_x$  
	\begin{align*}
	T_x^{n_\xi(x)}(B_x(z,\xi))=\jl_{\ta^{n_\xi(x)}(x)}.
	\end{align*}		
\end{itemize}
\begin{remark}\label{rem: cont inverse branch}
	Note that the measurably expanding condition implies that $T_x\rvert_{B_x(z,\eta_x)}$ is injective for each $(x,z)\in \cJ$. Furthermore, given that the spaces $\cJ_x$ are compact, we have that $\deg(T_x)$ is finite for each $x\in X$.
	Considering additionally the uniform openness condition we see that for every $(x,z)\in\cJ$ there exists a unique continuous inverse branch 
	\begin{align*}
		T_z^{-1}:B_{\ta(x)}(T_x(z),\xi)\lra B_x(z,\eta_x)
	\end{align*}
	of $T_x$ which sends $T_x(z)$ to $z$.
\end{remark}
The map $T$ is called \textit{uniformly expanding} if it satisfies the following additional properties:
\begin{enumerate}
	\item $\gamma_*:=\inf_{x\inX}\gamma_x>1$,
	\item $\deg(T):=\sup_{x\inX}\deg(T_x)<\infty$,
	\item $n_{\xi*}:=\sup_{x\inX}n_\xi(x)<\infty$.
\end{enumerate}
In addition to the various fiberwise and global Banach spaces defined in Section \ref{sec:potentials}, we will find the following definition useful. 

\begin{comment}
For each $x\inX$, we let $\sC(\jl_x)$ denote the space of all continuous functions $u_x:\jl_x\lra\RR$, and let $\sC(\jl)$ be the set of functions $u:\jl\lra\RR$ such that for $m$-a.e. $x\inX$ the functions $u_x:=u\rvert_{\jl_x}\in\sC(\jl_x).$  Let $\sC^0(\jl)$ and $\sC^1(\jl)$ be subspaces of $\sC(\jl)$ such that for all $u\in\sC^0(\jl)$ the function $x\longmapsto\norm{u_x}_\infty$ is measurable and $w\in\sC^1(\jl)$ implies that 
\begin{align*}
\norm{w}_1:=\int_X\norm{w_x}_\infty \,dm(x)<\infty.
\end{align*}
Let $\sC_*^\infty$ denote the space of all $\sB$--measurable functions $u\in\sC(\jl)$ such that 
\begin{align*}
	\sup_{x\inX}\norm{u_x}_\infty<\infty.
\end{align*}
Clearly $\sC^\infty_*$ becomes a Banach space when coupled together with the norm $\absval{\makebox[1ex]{\textbf{$\cdot$}}}_\infty$ given by
\begin{align*}
\absval{u}_\infty=\int_X\norm{u_x}_\infty\, dm(x), \quad u\in\sC^\infty_*.
\end{align*}

Fix $\al\in(0,1]$. The $\al$--variation of a function $u_x\in\sC(\jl_x)$ is given by
\begin{align*}
v_\al(u_x)=\inf_{\varrho_x(z_1,z_2)\leq \eta}\set{H_x:\absval{u_x(z_1)-u_x(z_2)}\leq H_x \varrho_x^\al(z_1,z_2)}.
\end{align*}
We say that a function $u\in\sC^1(\jl)$ is $\al$--H\"older if there is a measurable function $H:X\lra[1,\infty)$, $x\longmapsto H_x$,  such that $\log H\in L^1(m)$ and $v_\al(u_x)\leq H_x$ for $m$--a.e. $x\in X$. Let $\sH_\al(\jl)$ be the collection of all such functions. Clearly
\begin{align*}
\norm{\spot}_\al=\norm{\spot}_\infty+v_\al(\spot)
\end{align*}  
is a norm on $\sH_\al(\jl)$.

content...
\end{comment}
First, recall that $u\in\sH_\al(\cJ)$ provided $u_x\in\cH_\al(\cJ_x)$ and there is a measurable function $H:X\lra[1,\infty)$, $x\longmapsto H_x$,  such that $\log H\in L^1(m)$ and $v_{x,\al}(u_x)\leq H_x$ for $m$--a.e. $x\in X$.
For fixed $H$ we let $\sH_\al(\jl,H)$ be all the functions $g\in\sH_\al(\jl)$ such that
\begin{align*}
v_{x,\al}(g_x)\leq H_x
\end{align*}  
for $m$-a.e. $x\inX$ and we say that such a function is $(H,\al)$--H\"older continuous over $\cJ$. This allows us to write 
\begin{align*}
	\sH_\al(\jl)=\union_{H\geq 1}\sH_\al(\jl,H).
\end{align*}
Now for each $H$ and each $\phi\in\sH_\al(\jl,H)$ set 
\begin{align}
\cQ_x:=\cQ_x(H)=\sum_{j=1}^\infty H_{\ta^{-j}(x)}(\gamma^j_{\ta^{-j}(x)})^{-\al}.\label{def: Q_x}
\end{align}
The following lemma tells us this function is measurable and provides a necessary, though technical, bound. 
\begin{lemma}[Lemma 2.3 of \cite{mayer_distance_2011}]\label{lem: 2.3 derm}
	The function $x\longmapsto\cQ_x$ is measurable and $m$--a.e. finite. Moreover, for every $\phi\in\sH^\al(\jl,H)$, 
	\begin{align}
	\absval{S_n\phi_x(T^{-n}_y(w_1))-S_n\phi_x(T^{-n}_y(w_2))}\leq\cQ_{\ta^n(x)}\varrho_{\ta^n(x)}^\al(w_1,w_2)\label{lemma 2.3 S_n variation Q_x meas}	
	\end{align}
	for all $n\geq 1$, a.e. $x\inX$, every $z\in\jl_x$, and all $w_1,w_2\in B_{\ta^n(x)}(T^n_x(z),\xi)$.% where $y=(x,z)$.
\end{lemma}
\begin{comment}
Now let $\sH_\al^*(\jl)$ be the set of all functions $u$ such that $u\in\sC^\infty_*$ and \begin{align*}
\sup_{x\inX}v_\al(u_x)<\infty.
\end{align*} 
For $u\in\sH_\al^*(\cJ)$, we set 
\begin{align*}
V_\al(u)=\int_X v_\al(u_x)\, dm(x),
\end{align*}
and then define the norm $\absval{\makebox[1ex]{\textbf{$\cdot$}}}_\al$ on $\sH_\al^*(\cJ)$ by 
\begin{align*}
\absval{u}_\al=\absval{u}_\infty+V_\al(u).
\end{align*}

Then $\sH_\al^*(\jl)$ is a Banach space when considered with the norm $\absval{\makebox[1ex]{\textbf{$\cdot$}}}_\al$. 

content...
\end{comment}
For $\tilde{H}\geq 0$ we denote by $\sH_\al^*(\jl,\tilde{H})$ to be the space of all functions $\phi\in\sH_\al^*(\jl)$ such  that 
\begin{align*}
	\sup_{x\inX}v_{x,\al}(\phi_x)\leq H_x\leq \tilde{H}.	
\end{align*}
If $T$ is uniformly expanding then we are able to take $\cQ_x$ defined in \eqref{def: Q_x} as some uniform constant $\cQ_\phi$, depending only upon $\phi$ and $\tilde{H}$ and no longer on $x\in X$, given by 
\begin{align}
\cQ:=\cQ_\phi:=\tilde{H}\sum_{j=1}^\infty \gm^{-\al j}=\frac{\tilde{H}\gm^{-\al}}{1-\gm^{-\al}}.\label{eqn: Qtilde}
\end{align}
Thus, for uniformly expanding systems, we can rewrite Lemma \ref{lem: 2.3 derm} using the following lemma of \cite{mayer_distance_2011}.
\begin{lemma}[Lemma 3.31 of \cite{mayer_distance_2011}]\label{lem 3.31 DERM}
	For every $\phi\in\sH_\al^*(\jl,\tilde{H})$, 
	\begin{align}
	\absval{S_n\phi_x(T^{-n}_y(w_1))-S_n\phi_x(T^{-n}_y(w_2))}\leq\cQ_\phi\varrho_{\ta^n(x)}^\al(w_1,w_2)\label{ineq: uni exp var S_n Q }	
	\end{align}
	for all $n\geq 1$, a.e. $x\inX$, every $z\in\jl_x$, and all $w_1,w_2\in B_{\ta^n(x)}(T^n_x(z),\xi)$.
\end{lemma} 
Now we discuss the established thermodynamic formalism for distance expanding random systems.  We begin by defining the transfer operator. 
\begin{definition}
	Fix $\phi\in\sH_\al(\jl)$ and for each $x\in X$ define the operator\\ $\tr_x:=\tr_{\phi,x}:\cC(\jl_x)\lra\cC(\jl_{\ta(x)})$ by 
	\begin{align*}	
	\tr_x(u_x)(w)=\sum_{z\in T^{-1}_x(w)}u_x(z)e^{\phi_x(z)}.
	\end{align*}
	Clearly, $\tr_x$ is a positive bounded linear operator with norm bounded by 
	\begin{align*}
	\norm{\tr_x}_{x,\infty}\leq\deg(T_x)e^{\absval{\phi}_\infty}.
	\end{align*}
	Iterating the transfer operator, for each $n\in\NN$ we see that $\tr^n_x:\cC(\jl_x)\lra\cC(\jl_{\ta^n(x)})$ is given by 
	\begin{align*}
	\tr_x^n(u_x)(w)=\sum_{z\in T^{-n}_x(w)}u_x(z)e^{S_n\phi_x(z)}, \quad w\in\cJ_{\ta^n(x)}.
	\end{align*}
	We can also define the global operator $\tr:\sC(\jl)\lra\sC(\jl)$ by 
	\begin{align*}
	(\tr u)_x:=\tr_{\ta^{-1}(x)}u_{\ta^{-1}(x)}.
	\end{align*}
	Denote by $\tr_x^*$ the dual operator $\tr_x^*:\cC^*(\cJ_{\ta(x)})\lra\cC^*(\cJ_x)$, where $\cC^*(\cJ_x)$ denotes the dual space of $\cC(\cJ_x)$. 
	\begin{theorem}[Theorem 3.1 (1) of \cite{mayer_distance_2011}]
		There exists a unique family of probability measures $\nu_x\in\cP(\jl_x)$ such that 
		\begin{align*}
		\tr_x^*\nu_{\ta(x)}=\lm_x\nu_x \quad \text{ for } m-\text{a.e. }x\inX,
		\end{align*}
		where 
		\begin{align*}
		\lm_x:=\tr_x^*(\nu_{\ta(x)})(\ind)=\nu_{\ta(x)}(\tr_x\ind).
		\end{align*}
	\end{theorem}
	Define the normalized operator $\tr_{0,x}:\cC(\jl_x)\lra\cC(\jl_{\ta(x)})$ for each $x\in X$ by 
	\begin{align*}
	\tr_{0,x}=\lambda_x^{-1}\tr_x. 
	\end{align*}
	where $\lambda_x$ and $\nu_x$ are such that $\tr_x^*\nu_{\ta(x)}=\lambda_x\nu_x$, where $\rho_x$ is such that $\tr_{0,x}(\rho_x)=\rho_{\ta(x)}$. Let $\mu_x=\rho_x\nu_x$.
\end{definition}
For $p>2$ fix a function $g\in\sH_\al^p(\cJ)$. Then for each $r\in\RR$ and $x\in X$, define the perturbed operator $\tr_{r,x}:\cC(\jl_x)\lra\cC(\jl_{\ta(x)})$ by 
\begin{align*}
\tr_{r,x}u_x=\tr_{0,x}(e^{irg_x}\cdot u_x), 
\end{align*} 
and the global perturbed operator by 
\begin{align*}
\left(\tr_r u\right)_x=\tr_{r,\ta^{-1}(x)}u_{\ta^{-1}(x)}.		
\end{align*}
\begin{remark}
	Note that by our assumption that the function $g:\jl\lra\RR$ is bounded and H\"{o}lder we have that the process $g\circ T$ is in $L^p$ for any $p>2$.
\end{remark}
Now let $r_{n-1},\dots,r_0\in\RR$. Then 
\begin{align*}
\int_{\cJ_{\ta^n(x)}} \tr_{r_{n-1},\ta^{n-1}(x)}\circ\cdots\circ\tr_{r_0,x}( \rho_x)\,d\nu_{\ta^{n}(x)}
&=\int_{\cJ_{\ta^n(x)}} \tr^n_{0,x}(e^{i\sum_{j=0}^{n-1}r_jg_{\ta^j(x)}\circ T^j_x}\cdot  \rho_x)\,d\nu_{\ta^{n}(x)}\\
&=\int_{\cJ_x} e^{i\sum_{j=0}^{n-1}r_jg_{\ta^j(x)}\circ T^j_x}\cdot  \rho_x \,d\nu_x\\
&=\int_{\cJ_x} e^{i\sum_{j=0}^{n-1}r_jg_{\ta^j(x)}\circ T^j_x}\, d\mu_x
\end{align*}

Integrating with respect to $\nu$ then gives 
\begin{align*}
\nu\left(\tr_{r_{n-1}}\circ\cdots\circ\tr_{r_0}( \rho)\right)
&=\int_X\int_{\jl_{\ta^n(x)}}\tr_{r_{n-1},\ta^{n-1}(x)}\circ\cdots\circ\tr_{r_0,x}( \rho_x)\,d\nu_{\ta^{n}(x)}\,dm(x)\\
&=\int_X\int_{\jl_x} e^{i\sum_{j=0}^{n-1}r_jg_{\ta^j(x)}\circ T^j_x}\,d\mu_x\,dm(x)\\
&=\cE_\mu\left(e^{i\sum_{j=0}^{n-1}r_jg\circ T^j}\right).	
\end{align*}

For fixed $r\in\RR$ and any $n\in\NN$ we see that we can write the iterates of $\tr_r$ as  
\begin{align*}
\tr_{r,x}^n(u_x)=\tr_{0,x}^n(e^{irS_ng_x}\cdot u_x).
\end{align*}

Now we endeavor to show that for sufficiently small values of $\absval{r}$, the operator $\tr_r$ is bounded in the $\al$--norm and the sup norm. We start by applying the following lemma. 
\begin{lemma}[Lemma 3.8 of \cite{mayer_distance_2011}]\label{lem 3.8 DERM}
	For all $w_1,w_2\in \jl_x$ and $n\geq 1$
	\begin{align*}
	\frac{\tr_{0,\ta^{-n}(x)}^n\ind(w_1)}{\tr_{0,\ta^{-n}(x)}^n\ind(w_2)}=\frac{\tr_{\ta^{-n}(x)}^n\ind(w_1)}{\tr_{\ta^{-n}(x)}^n\ind(w_2)}\leq C_\phi(x),
	\end{align*}
	where $C_\phi$ is given by
	\begin{align*}
	C_\phi(x):=e^{\cQ_{\ta^{-j}(x)}}\deg(T^j_{\ta^{-j}(x)})\max\set{\exp(2\norm{S_k\phi_{\ta^{-k}(x)}}_{x,\infty}):0\leq k\leq j}\geq 1.
	\end{align*}
	If in addition we have that $\varrho_x(w_1,w_2)\leq \xi$, then
	\begin{align*}
	\frac{\tr_{0,\ta^{-n}(x)}^n\ind(w_1)}{\tr_{0,\ta^{-n}(x)}^n\ind(w_2)}\leq \exp(\cQ_x\varrho_x^\al(w_1,w_2)).
	\end{align*}
	Moreover, 
	\begin{align*}
	\frac{1}{C_\phi(x)}\leq \tr_{0,\ta^{-n}(x)}^n\ind(w)\leq C_\phi(x)
	\end{align*}
	for every $w\in\jl_x$ and $n\geq 1$.
\end{lemma}

However, we note that for uniformly expanding systems, in view of \eqref{eqn: Qtilde}, there is some $C>1$ such that for every $x\inX$ and $n\in\NN$
\begin{align}
C^{-1}\leq\tr_{0,x}^n\ind_x\leq C.\label{eqn:L0 bounds}
\end{align}
In light of \eqref{eqn:L0 bounds}, then for $r\in\RR$ and $u\in\sH_\al^*(\jl)$ we can estimate 
\begin{align}
\norm{\tr_{r,x}^n(u_x)}_{\ta^n(x),\infty}&=\norm{\tr_{0,x}^n(e^{irS_ng}\cdot u_x)}_{\ta^n(x),\infty}
\leq\norm{\tr_{0,x}^n(\absval{u_x})}_{\ta^n(x),\infty}\nonumber\\
&\leq\norm{u_x}_{x,\infty}\cdot\norm{\tr^n_{0,x}(\ind_x)}_{\ta^n(x),\infty}
\leq C\norm{u_x}_{x,\infty}. \label{eqn:infty bound on Lr}
\end{align}
\begin{comment}
Integrating over $X$ gives 
\begin{align}
\absval{\tr_r^nu}_\infty\leq C\absval{u}_\infty.\label{ineq:Lrinfty bdd}
\end{align}

content...
\end{comment}
Now to show that $\norm{\tr^n_{r,x}u_x}_{\ta^n(x),\al}\leq C\norm{u_x}_{x,\al}$ we must first show that this inequality holds on positive cones of H\"older functions. For $s\geq 1$ and $x\inX$ let 
\begin{align*}
\Lambda^s_x=\{h_x\in\cC(\jl_x):h_x\geq 0,\nu_x(h_x)=1, \text{ and } &h_x(w_1)\leq \exp(s\cQ_x\varrho_x^\al(w_1,w_2))\cdot h_x(w_2)\\
&\text{for all }w_1,w_2\in\jl_x \text{ with }\varrho_x(w_1,w_2)\leq\xi\}.
\end{align*}
The first of the following two lemmas shows that $\Lambda_x^s\sub\cH_\al(\jl_x)$, and moreover, provides a bound on the variation of such functions, while the second lemma is a sort of converse to the first.
\begin{lemma}[Lemma 3.11 of \cite{mayer_distance_2011}]\label{lem 3.11 DERM}
	If $h_x\geq 0$ and for all $w_1,w_2\in\jl_x$ with $\varrho_x(w_1,w_2)\leq \xi$ we have 
	\begin{align*}
	h_x(w_1)\leq e^{s\cQ_x\varrho^\al(w_1,w_2)}h_x(w_2),
	\end{align*}
	then 
	\begin{align*}
	v_{x,\al}(h_x)\leq s\cQ_x(\exp(s\cQ_x\xi^\al))\xi^\al\norm{h_x}_{x,\infty}.
	\end{align*}
\end{lemma}

\begin{lemma}[Lemma 3.13 of \cite{mayer_distance_2011}]\label{lem 3.13 DERM}
	If $u_x\in\cH_\al(\jl_x)$ and $u_x\geq 0$, then the function defined by
	\begin{align*}
	h_x:=\frac{u_x+v_{x,\al}(u_x)/\cQ_x}{\nu_x(u_x)+v_{x,\al}(u_x)/\cQ_x}\in\Lambda_x^1.
	\end{align*}
\end{lemma}
The next lemma of \cite{mayer_distance_2011} establishes the invariance of the cones $\Lm_x^s$ with respect to the family of normalized operators $\tr_{0,x}$. 
\begin{lemma}[Lemma 3.14 of \cite{mayer_distance_2011}]\label{lem 3.14 DERM}
	Let $h_x\in\Lm_x^s$. Then for every $n\geq 1$ and $w_1,w_2\in\jl_x$ with $\varrho(w_1,w_2)\leq\xi$ we have
	\begin{align}
	\frac{\tr_{0,x}^nh_x(w_1)}{\tr_{0,x}^nh_x(w_2)}\leq\exp(s\cQ_{\ta^n(x)}\varrho^\al(w_1,w_2)). 
	\end{align}
	Consequently, $\tr_{0,x}^n(\Lm_x^s)\sub\Lm_{\ta^n(x)}^s$ for a.e. $x\inX$ and all $n\geq 1$.
\end{lemma}
Now we wish to estimate the values $\norm{\tr_{r,x}^nu_x}_{\ta^n(x),\al}$.
\begin{lemma}
There exists $C\geq1$ such that for all $n\in\NN$, all $r\in\RR$, $m$-a.e. $x\inX$, and all $u\in\sH_\al(\cJ)$ we have
\begin{align*}
%\norm{\tr_{r,x}^nf_x}_{\ta^n(x),\infty}\leq C\norm{f_x}_{x,\infty} \spand
 \norm{\tr_{r,x}^nu_x}_{\ta^n(x),\al}\leq C\norm{u_x}_{x,\al}.
\end{align*}		
\end{lemma}
\begin{proof}
		
In comparing $u_x$ and $h_x$ we note that 
\begin{align}
v_{x,\al}(h_x)&=\frac{\cQ_x}{\cQ_x\nu_x(u_x)+v_{x,\al}(u_x)}v_{x,\al}(u_x), \quad \text{and}\quad \nonumber\\
\norm{h_x}_{x,\infty}&=\absval{\frac{\cQ_x}{\cQ_x\nu_x(u_x)+v_{x,\al}(u_x)}}\norm{u_x}_{x,\infty}+\absval{\frac{v_{x,\al}(u_x)}{\cQ_x\nu_x(u_x)+v_{x,\al}(u_x)}}.\label{eqn:hx infty norm}.
\end{align} 
Now for $u_x\in\cH_\al(\jl_x)$ we can write $u_x=u_x^+-u_x^-$ where $u_x^+,u_x^-$ are both in $\cH_\al(\jl_x)$ and nonnegative. Letting 
\begin{align*}
h_x^+=\frac{u_x^++v_{x,\al}(u_x^+)/\cQ_x}{\nu_x(u_x^+)+v_{x,\al}(u_x^+)/\cQ_x}\quad \text{ and } \quad 
h_x^-=\frac{u_x^-+v_{x,\al}(u_x^-)/\cQ_x}{\nu_x(u_x^-)+v_{x,\al}(u_x^-)/\cQ_x},
\end{align*}then for $h_x^+$ we have the following estimate
\begin{align*}
v_{x,\al}(\tr_{0,x}^n(u_x^+))&= \frac{\cQ_x\nu_x(u_x^+)+v_{x,\al}(u_x^+)}{\cQ_x}\left(v_{x,\al}(\tr_{0,x}^nh_x^+)-\frac{v_{x,\al}(u_x^+)}{\cQ_x\nu_x(u_x^+)+v_{x,\al}(u_x^+)}v_{x,\al}(\tr^n_{0,x}\ind_x)\right)\\ 
&\leq 
\cQ_x(\exp(\cQ_x\xi^\al))\xi^\al\frac{\cQ_x\nu_x(u_x^+)+v_{x,\al}(u_x^+)}{\cQ_x}
\left(\norm{\tr_{0,x}^nh_x^+}_{\ta^n(x),\infty}+\norm{\tr^n_{0,x}\ind_x}_{\ta^n(x),\infty}\right)\\
&\leq\xi^\al\exp(\cQ_x\xi^\al)\left(\cQ_x\nu_x(u_x^+)+v_{x,\al}(u_x^+)\right)\left(\norm{\tr_{0,x}^nh_x^+}_{\ta^n(x),\infty}+C\right)\\
&\leq\xi^\al\exp(\cQ_x\xi^\al)\left(\cQ_x\nu_x(u_x^+)+v_{x,\al}(u_x^+)\right)\left(\frac{\cQ_x}{\cQ_x\nu_x(u_x^+)+v_{x,\al}(u_x^+)}\norm{\tr_{0,x}^nu_x^+}_{\ta^n(x),\infty}+C\right)\\
&\leq\max(\cQ_x,1)\xi^\al\cQ_x\exp(\cQ_x\xi^\al)\norm{\tr_{0,x}^nu_x^+}_{\ta^n(x),\infty}+C\norm{u_x^+}_{x,\al}.
\end{align*}
We get a similar bound for $v_{\ta^n(x),\al}(\tr_{0,x}^nu_x^-)$. Now if $u_x\in\sH_\al^*(\cJ)$ then so are $u_x^+,u_x^-$ and furthermore $\norm{u_x^+}_{x,\infty}, \norm{u_x^-}_{x,\infty}\leq\norm{u_x}_{x,\infty}$ and $v_{x,\al}(u_x^+),v_{x,\al}(u_x^-)\leq v_{x,\al}(u_x)$. If $T$ is uniformly expanding then, using \eqref{eqn: Qtilde}, we are able to simplify the above estimate to  
\begin{align*}
v_{\ta^n(x),\al}(\tr_{0,x}^n(u_x^+))&\leq \max(\cQ_u,1)\xi^\al\cQ_u\exp(\cQ_u\xi^\al)\norm{\tr_{0,x}^nu_x^+}_{\ta^n(x),\infty}+M\norm{u_x^+}_{x,\al}\\
&\leq C\left(\norm{\tr_{0,x}^nu_x^+}_{\ta^n(x),\infty}+\norm{u_x^+}_{x,\al}\right).
\end{align*}
Again, we obtain a similar estimate for $v_{\ta^n(x),\al}(\tr_{0,x}^nu_x^-)$. Combining \eqref{eqn:hx infty norm} and \eqref{eqn:L0 bounds}, we see
\begin{align}
v_{\ta^n(x),\al}(\tr_{0,x}^n(u_x))&\leq v_{\ta^n(x),\al}(\tr_{0,x}^nu_x^+)+v_{\ta^n(x),\al}(\tr_{0,x}^nu_x^-)\nonumber\\
&\leq C\left(\norm{u_x^+}_{x,\infty}+\norm{u_x^+}_{x,\al}+\norm{u_x^-}_{x,\infty}+\norm{u_x^-}_{x,\al}\right)\nonumber\\ 
&\leq C\norm{u_x}_{x,\al}.\label{eqn: valpha bound}
\end{align}
Now we are ready to calculate $\norm{\tr^n_{r,x}u_x}_{\ta^n(x),\al}$. Combining the bounds in \eqref{eqn:infty bound on Lr} and \eqref{eqn: valpha bound} we see
\begin{align*}
\norm{\tr^n_{r,x}u_x}_{\ta^n(x),\al}=\norm{\tr_{0,x}^n(e^{irS_ng_x}\cdot u_x)}_{\ta^n(x),\al}\leq C\norm{u_xe^{irS_ng_x}}_{x,\al}. 
\end{align*}
As $\sH_\al^*(\cJ)$ is a Banach algebra we can write the last inequality as 
\begin{align*}
\norm{\tr^n_{r,x}u_x}_{\ta^n(x),\al}&\leq C\norm{u_x}_{x,\al}\norm{e^{irS_ng_x}}_{x,\al}.
\end{align*}
It suffices to estimate $\norm{e^{irS_ng_x}}_{x,\al}$, but as $\norm{e^{irS_ng_x}}_{x,\infty}=1$, we need only estimate $v_{x,\al}(e^{irS_ng_x})$. Since $\absval{e^{iz_1}-e^{iz_2}}\leq\absval{z_1-z_2}$, 
for all $z_1,z_2\in\RR$,  we see that $v_{x,\al}(e^{irS_ng_x})\leq v_{x,\al}(rS_ng_x)$. Thus it now suffices to estimate $v_{x,\al}(rS_ng_x)$. Lemma \ref{lem 3.31 DERM}, shows that
\begin{align*}
v_{x,\al}(e^{irS_ng_x})\leq v_{x,\al}(rS_ng_x)\leq \absval{r}\cQ.
\end{align*} 
In particular, we have the estimate 
\begin{align*}
\norm{\tr^n_{r,x}u_x}_{\ta^n(x),\al}&\leq C\norm{u_x}_{x,\al}(1+\absval{r}\cQ), 
\end{align*}
which finishes the proof.
\end{proof}
\begin{comment}
Integrating over $X$ and considering only $r\in\RR$ with $\absval{r}\leq 1$, gives the desired inequality 
\begin{align*}
\absval{\tr_r^nu}_\al\leq C\absval{u}_\al, \quad n\in\NN.
\end{align*}

content...
\end{comment}
Now we define the fiberwise operator $Q_x:\cC(\jl_x)\lra\cC(\jl_{\ta(x)})$ by 
\begin{align*}
Q_x(u_x)=\int_{\cJ_x}u_x\,d\nu_x\cdot \rho_{\ta(x)}.
\end{align*}
For each $n\in\NN$ we can write the iterates $Q_x^n:\cC(\cJ_x)\lra\cC(\cJ_{\ta^n(x)})$ of $Q_x$ as 
\begin{align*}
Q_x^n(u_x)=\int_{\cJ_x}u_x\,d\nu_x\cdot \rho_{\ta^n(x)}.
\end{align*}
Thus we see that for each $k\in\NN$ and $x\inX$
\begin{align*}
\norm{\left((\tr_0^k-Q^k)(u)\right)_x}_{x,\infty}&=\norm{\left(\tr_{0,\ta^{-k}(x)}^k-Q^k_{\ta^{-k}(x)}\right)(u_{\ta^{-k}(x)})}_{x,\infty}\\
&=\norm{\tr_{0,\ta^{-k}(x)}^k(u_{\ta^{-k}(x)})-\int_{\cJ_{\ta^{-k}(x)}}u_{\ta^{-k}(x)} \,d\mu_{\ta^{-k}(x)}\cdot\ind_x  }_{x,\infty}.
\end{align*}
By the proof of Lemma 3.18 of \cite{mayer_distance_2011}, in the case of uniformly expanding random maps, we see that
\begin{align*}
\norm{\left((\tr_0^k-Q^k)(u)\right)_x}_{x,\infty}&\leq \left(\nu_x(u_x)+2\frac{v_{x,\al}(u_x)}{Q_x}\right)C\kp^k \\
&\leq\left(\norm{u_x}_{x,\infty}+2v_{x,\al}(u_x)\right)C\kp^k\\
&\leq 2\norm{u_x}_{x,\al} C\kp^k
\end{align*}
for some positive constant $\kp<1$. 

Thus we have shown that for uniformly expanding systems conditions \eqref{assum 1}--\eqref{assum 4} of Theorem \ref{thm: modified gouezel asip} hold. In order for the remaining two conditions to hold we will need to require more structure. In the next section we discuss a class of systems, first described by Denker and Gordin in \cite{denker_gibbs_1999}, which fit within the framework of the uniformly expanding systems which we have just discussed. In the same manner as in \cite{mayer_distance_2011}, we shall call refer to these systems as DG and DG*--systems. 
%By our assumption that the function $g:\jl\lra\RR$ is bounded and H\"{o}lder we have that the process $g\circ T$ is in $L^p$ for any $p>2$.

\section{DG*--Systems}

\begin{comment}
\bc\textbf{THINGS TO FIX}\ec
\begin{itemize}
	\item fix all constants and keep track don't just use C
	\item figure out which exponent is larger alpha or beta. if neither use alpha' is the max or something
	\item fix all exponents and keep track from the first things proved. many later alpha may need to be betas. (or whichever one is larger)
	
\end{itemize}

--------------------------------------------------

content...
\end{comment}

In \cite{denker_gibbs_1999} Denker and Gordin first established the existence and uniqueness of conformal Gibbs measures for DG--systems. Then in \cite{mayer_distance_2011} Mayer, Skorulski, and Urba\'nski were able to cast these systems as uniformly expanding random systems, which they called DG*--systems, meaning that the full thermodynamic formalism they developed there applies to these DG*--systems. In particular, we see that the results of Section \ref{sec: RDEM} apply, and in particular, the spectral gap property holds for these systems. Thus we have only to check conditions \eqref{assum 5} and \eqref{assum 6} of Theorem \ref{thm: modified gouezel asip} hold. However, condition \eqref{assum 5} has been shown to hold.
In this section we introduce these systems and show that an ASIP holds for such systems. We begin with a definition.

\begin{definition}\label{def: DG system}
	Suppose $(X_0, d_{X_0})$ and $(Z_0,d_{Z_0})$ are compact metric spaces and that $\ta_0:X_0\lra X_0$ and $T_0:Z_0\lra Z_0$ are open topologically exact distance expanding mappings in the sense of \cite{przytycki_conformal_2010}. Assume that $T_0$ is a skew product over $Z_0$, that is, for each $x\in X_0$ there exists a compact metric space $\cJ_x$ such that $Z_0=\union_{x\in X_0}\set{x}\times\cJ_x$. Further assume that the map $\ta_0$ is Lipschitz, i.e. there exists $L>0$ such that 
	\begin{align*}
		d_{X_0}(\ta_0(x),\ta_0(x'))\leq L d_{X_0}(x,x')
	\end{align*}
	for all $x,x'\in X_0$, and that there exists $\xi,\xi_1>0$ such that for all $x,x'\in X_0$ with $d_{X_0}(x,x')<\xi_1$ there exist $y\in\cJ_x$ and $y'\in\cJ_{x'}$ such that 
	\begin{align*}
		d_{Z_0}((x,y),(x',y'))<\xi.
	\end{align*}
	Finally we assume that the projection $\pi:Z_0\lra X_0$ onto the first coordinate given by 
	\begin{align*}
		\pi(x,y)=x
	\end{align*}
	is an open mapping and that the following diagram commutes:
	\bc
\begin{tikzcd}
	Z_0 \arrow[r, "T_0"] \arrow[d, "\pi"]
	& Z_0 \arrow[d, "\pi" ] \\
	X_0 \arrow[r, "\ta_0" ]
	& X_0
\end{tikzcd}
	\ec
	We call the system $(T_0,Z_0,\ta_0,X_0)$ a DG--system. 	
\end{definition}
In what follows we will assume that there is some metric space $(Y,d_Y)$ such that $\cJ_x\sub Y$ for each $x\in X$. In this case we see that $Z_0\sub X\times Y$ and we may take $d_{Z_0}$ to be the natural product metric given by 
\begin{align*}
	d_{Z_0}((x,y),(x',y'))=d_{X_0}(x,x')+d_Y(y,y').
\end{align*}
\begin{remark}
	Three things to notice about the definition of DG--systems presented above.
	\begin{itemize}
		\item For each $x\in X_0$ we have that 
		\begin{align*}
			T_0(\set{x}\times\cJ_x)\sub\set{\ta_0(x)}\times\cJ_{\ta_0(x)},
		\end{align*}
		giving rise to the map $T_x:\cJ_x\lra\cJ_{\ta_0(x)}$.
		\item Since $T_0$ is distance expanding in the sense of \cite{przytycki_conformal_2010} the conditions of uniform openness, measurably expanding, measurability of the degree, and topological exactness from the definition of random distance expanding mappings in Section \eqref{sec: RDEM} all hold for constants $\gm_x\geq \gm>1$, $\deg(T_x)\leq N_1<\infty$, and $n_\xi=n_\xi(x)$ independent of $x$. 
		\item The function $\ta_0:X_0\lra X_0$ need not be invertible, meaning that we are not quite able to apply the theory of uniformly expanding random mappings which we described earlier. 
	\end{itemize}
\end{remark}
In order to rectify the complications involving $\ta_0$ we define DG*--systems by turning to Rokhlin's natural extension, i.e. the projective limit, $\ta:X\lra X$ of $\ta_0:X_0\lra X_0$.
\begin{definition}
	Assume that we are given a DG--system $(T_0,Z_0,\ta_0,X_0)$ as defined above. We further assume that the space $X_0$ comes coupled with a Borel probability $\ta_0$--invariant ergodic measure $m_0$ and a H\"older continuous potential $\phi:Z_0\lra\RR$. Define the space 
	\begin{align*}
		X=\set{(x_n)_{n\leq 0}:\ta_0(x_n)=x_{n+1} \text{ for all }n\leq -1}
	\end{align*}
	and the map $\ta:X\lra X$ by 
	\begin{align*}
		\ta((x_n)_{n\leq 0})=(\ta_0(x_n))_{n\leq 0}.
	\end{align*}
	In this case we have that $\ta:X\lra X$ is invertible and that the diagram
	\bc
	\begin{tikzcd}
		X \arrow[r, "\ta"] \arrow[d, "p"]
		& X \arrow[d, "p" ] \\
		X_0 \arrow[r, "\ta_0" ]
		& X_0
	\end{tikzcd}
	\ec
	commutes where $p:X\lra X_0$ is defined by 
	\begin{align*}
		p((x_n)_{n\leq 0})=x_0.
	\end{align*}
	Now since $m_0$ is a $\ta_0$--invariant ergodic measure, there exists a unique $\ta$--invariant probability measure $m$ on $X$ such that $m\circ \pi^{-1}=m_0$.  Define the set
	\begin{align*}
		\cJ=\union_{x\in X}\set{x}\times \cJ_{x_0}
	\end{align*}
	and the map $T:\cJ\lra \cJ$ by 
	\begin{align*}
		T(x,y)=(\ta(x),T_{x_0}(y)).
	\end{align*}
	The metrics $d_{X_0}$ and $d_{Z_0}$ extend naturally to metrics $d_X$ and $d_\cJ$ on $X$ and $\cJ$, which are defined by 
	\begin{align*}
		d_X(x,x')=\sum_{n=0}^\infty 2^{-n}d_{X_0}(x_{-n},x'_{-n}), \quad x,x'\in X
	\end{align*}
	and 
	\begin{align*}
		d_\cJ((x,y),(x',y'))=d_Y(y,y')+d_{X}(x,x'), \quad (x,y),(x',y')\in\cJ
	\end{align*}
	respectively.  We let $B_\cJ((x,w),r)$ denote the ball of radius $r>0$ centered at the point $(x,w)\in\cJ$ with respect to the metric $d_\cJ$. In what follows we will assume that $m=m_\psi$ is a $\ta$--invariant Gibbs measure for some continuous H\"older potential $\psi$ on $X$, having nothing to do with our (fiberwise) potential $\phi$ or density function $\rho$. The system $(T,\cJ,m,\ta,X)$ is then called a DG*--system.
\end{definition}
\begin{remark}\label{rem:DG* system ineq}
Note that because DG*--systems are uniformly expanding random systems, we have that there exists $C>1$ such that 
\begin{itemize}
	\item $C^{-1}\leq \rho_x\leq C$ for all $x\inX$,
	\item $C^{-1}\leq \tr_{0,\ta^{-n}(x)}^n\ind_{\ta^{-n}(x)}\leq C$ for all $x\in X$ and $n\in\NN$. 
\end{itemize}
\end{remark}

	Notice that since the base dynamical system $(X,\ta,m)$ is distance expanding in the sense of \cite{przytycki_conformal_2010} and that $m$ is an invariant Gibbs measure, the base system exhibits an exponential decay of correlations for functions $F$ and $G$ so long as one of the two is H\"older continuous on $X$ and the other is integrable with respect to $m$.  More precisely, we have the following, which is a consequence of Theorem 5.4.9 of \cite{przytycki_conformal_2010}.
	\begin{theorem}
		Let $(T,\cJ,m,\ta,X)$ be a DG*--system such that $m$ is an invariant Gibbs measure. Then there exists $C\geq1$ and $\kp<1$ such that for all $F\in\HH_\bt$ and all $G\in L^1(m)$ we have
		\begin{align*}
			\absval{\int_X (G\circ \ta^{-n})\cdot F\, dm-\int_X G\, dm\cdot \int_XF\,dm}\leq C\kp^n\norm{F}_\HH \norm{G}_{L^1(m)}
		\end{align*}
	\end{theorem}
	In other words, condition \eqref{assum 5} of Theorem \ref{thm: modified gouezel asip} is satisfied for DG*--systems. 

The following theorem addresses condition \eqref{assum 6} of Theorem \ref{thm: modified gouezel asip} and is due to Denker and Gordin. 
\begin{theorem}[Theorem 2.10 of \cite{denker_gibbs_1999}]\label{thm: DG Holder contin}
	Let $(T_0,Z_0,\ta_0,X_0)$ be a DG--system with random Gibbs measures $(\nu_x)_{x\in X_0}$. 
	If $f:Z_0\lra\RR$ is $(D,\al)$--H\"older continuous, for any $\al\in(0,1]$ and $D>0$,  then there exists $\bt_\al\in(0,1]$ and $C_\al>0$ such that the function
	\begin{align*}
		x\longmapsto\int_{\cJ_x}f_x(z)\,d\nu_x(z)
	\end{align*} is $(C_\al D,\bt_\al)$--H\"older continuous on $X_0$. In particular the function 
	\begin{align*}
		x\longmapsto\lm_x
	\end{align*} 
	is $\bt_\al$--H\"older continuous on $X_0$.
\end{theorem}
\begin{remark}
	Note that the same theorem applies for DG*--systems. Theorem 8.12 of \cite{mayer_distance_2011} reproves the special case of the previous theorem for the function $x\longmapsto\lm_x$ in the case of DG*--systems. 
\end{remark}
We shall now prove the following. 
\begin{proposition}\label{lem: perturbed int holder}
	Suppose that $u\in\sH_\tau(\cJ)$, for $\tau\in(0,1]$, and that there is $C>1$ such that 
	\begin{align*}
		C^{-1}\leq u_x\leq C 
	\end{align*}
	for all $x\inX$. Given $g\in \sH_\al^p(\cJ)$ and $\phi\in\sH_\al(\cJ)$, define the transfer operators as before. Then there exists $\bt\in(0,1]$, depending only on $\al$ and $\tau$, such that for each $n\in\NN$ and each $r_0,\dots,r_{n-1}\in\RR$ with each $\absval{r_j}<\ep_0$, for some $\ep_0>0$, the function
	\begin{align*}
		x\longmapsto\int_{\cJ_x}\tr_{r_{n-1},\ta^{-1}(x)}\circ\cdots\circ\tr_{r_0,\ta^{-n}(x)}(u_{\ta^{-n}(x)})\,d\nu_x
	\end{align*}  
	is in $\HH_\bt$.	
	Moreover, we have that 
	\begin{align*}
		\norm{\int_{\cJ_x}\tr_{r_{n-1},\ta^{-1}(x)}\circ\cdots\circ\tr_{r_0,\ta^{-n}(x)}(u_{\ta^{-n}(x)})\,d\nu_x}_\HH\leq C
	\end{align*} 
	independent of the choice of $n$ or the $r_j$.
\end{proposition}

\begin{proof}[Proof of Proposition \ref{lem: perturbed int holder}]	

Fix $n\in\NN$ and let $(x,w)\in\cJ$.  Then for each $z\in T^{-n}_x(w)$ there is a unique continuous inverse branch  
\begin{align*}
	T_{(\ta^{-n}(x),z)}^{-n}:B_{\cJ}((x,w),\xi)\longrightarrow B_{\cJ}((\ta^{-n}(x),z),\xi)
\end{align*}
which sends the point $(x,w)$ to $(\ta^{-n}(x),z)$. Similarly, for $(x',w')\in B_{\cJ}((x,w),\xi)$, there is a unique continuous inverse branch 
\begin{align*}
T_{(\ta^{-n}(x'),z`)}^{-n}:B_{\cJ}((x,w),\xi)\longrightarrow B_{\cJ}((\ta^{-n}(x),z),\xi)
\end{align*}
which sends the point $(x',w')$ to $(\ta^{-n}(x'),z')\in B_{\cJ}((\ta^{-n}(x),z),\xi)$. Thus for each $z\in T_x^{-n}(w)$ there is a unique and bijectively defined $z'\in T^{-n}_{x'}(w')$ such that $z'\in B_{\cJ}((\ta^{-n}(x),z),\xi)$.
Consequently, we may re--index the sum in the definition of the transfer operator allowing us to write
\begin{align}\label{frm:tran index}
\tr_{x'}^n(u_{\ta^{-n}(x')})(w')&=\sum_{z'\in T_{x'}^{-n}(w')}e^{S_n\phi_{\ta^{-n}(x')}(z')}u_{\ta^{-n}(x')}(z')\nonumber\\
&=\sum_{z\in T_{x}^{-n}(w)}e^{S_n\phi_{\ta^{-n}(x')}(z')}u_{\ta^{-n}(x')}(z').
\end{align} 
For $\zt\in (0,1]$ let $\bt_\zt$ be the number coming from Theorem \ref{thm: DG Holder contin}.
In order to prove Proposition \ref{lem: perturbed int holder} we wish to employ Theorem \ref{thm: DG Holder contin}, thus it suffices to show that 
\begin{align*}
	\tr_{r_{n-1}}\circ\cdots\circ\tr_{r_0}(u)\in\sH_\zt(\cJ).
\end{align*}
for some $\zt\in(0,1]$. To that end, we will first prove two lemmas, however we begin with an observation concerning the following definition. For each $n\in\NN$ and $r_0,\dots,r_{n-1}\in\RR$, set 
\begin{align*}
\ol{S}_nh_x(z):=\ol{S}_{n,r_0,\dots,r_{n-1}}h_x(z):=\sum_{j=0}^{n-1}r_j\cdot h_{\ta^j(x)}\circ T_x^j(z)
\end{align*}
for $h:\cJ\lra\RR$. Now for $h\in\sH_\zt(\cJ)$, for any $\zt\in(0,1]$, we have that
\begin{align*}
\absval{h(T^k(\ta^{-n}(x),z))-h(T^k(\ta^{-n}(x'),z'))}&\leq Cd_{\cJ}^\zt(T^k(\ta^{-n}(x),z),T^k(\ta^{-n}(x'),z'))\\
&\leq Cd_{\cJ}^\zt((x,w),(x',w'))\cdot \gm^{-\zt(n-k)}.
\end{align*}
Hence we have
\begin{align}
\absval{S_nh_{\ta^{-n}(x)}(z)-S_nh_{\ta^{-n}(x')}(z')}
&\leq Cd_{\cJ}^\zt((x,w),(x',w'))\cdot\sum_{k=0}^{n-1}\gm^{-\zt(n-k)}\nonumber\\
&\leq \frac{C}{1-\gm^{-\zt}}\cdot d_{\cJ}^\zt((x,w),(x',w')),\label{ineq: Sn DG*}
\end{align} 
and similarly we have 
\begin{align}
\absval{\ol{S}_nh_{\ta^{-n}(x)}(z)-\ol{S}_nh_{\ta^{-n}(x')}(z')}
%	&\leq Cd_{\cJ}^\al((x,w),(x',w'))\cdot\sum_{k=0}^{n-1}\gm^{-\al(n-k)}\\
&\leq \frac{C}{1-\gm^{-\zt}}\cdot d_{\cJ}^\zt((x,w),(x',w')).\label{ineq: Snbar DG*}
\end{align}
We now wish to prove the first of two lemmas.
\begin{lemma}\label{lem: tr op HC}
		The function $\tr^n\ind$ is $\al$--H\"older on $\cJ$.
\end{lemma}
\begin{proof}
		To see this we calculate 
		\begin{align*}
		&\absval{\tr_{\ta^{-n}(x)}^n\ind_{\ta^{-n}(x)}(w)-\tr_{\ta^{-n}(x')}^n\ind_{\ta^{-n}(x')}(w')}
		=\absval{\sum_{z\in T_{x}^{-n}}e^{S_n\phi_{\ta^{-n}(x)}(z)} - \sum_{z'\in T_{x'}^{-n}}e^{S_n\phi_{\ta^{-n}(x')}(z')} }\\
		&\leq \sum_{z\in T_{x}^{-n}(w)}\absval{e^{S_n\phi_{\ta^{-n}(x)}(z)} - e^{S_n\phi_{\ta^{-n}(x')}(z')}}\\
		&= \sum_{z\in T_{x}^{-n}(w)}\left( \absval{e^{S_n\phi_{\ta^{-n}(x')}(z')}}\absval{e^{S_n\phi_{\ta^{-n}(x)}(z) - S_n\phi_{\ta^{-n}(x')}(z')}-1} \right)\\
		&\leq C \cdot\sum_{z\in T_{x}^{-n}(w)}\left( \absval{e^{S_n\phi_{\ta^{-n}(x')}(z')}}\absval{S_n\phi_{\ta^{-n}(x)}(z) - S_n\phi_{\ta^{-n}(x')}(z')} \right)\\
		&\leq C\cdot\tr_{\ta^{-n}(x')}^n\ind_{\ta^{-n}(x')}(w')\cdot d_{\cJ}^\al((x,w),(x',w'))\\
		&\leq C\cdot d_{\cJ}^\al((x,w),(x',w')).
		\end{align*}
		This finishes the proof. 
\end{proof}
		Using the previous lemma we can now show the same for the normalized operator.
\begin{lemma}\label{lem: norm tr op HC}
		The function $\tr_0^n\ind$	is $\vkp$--H\"older on $\cJ$, where $\vkp:=\min\set{\al,\bt_\al}$.	
\end{lemma}
\begin{proof}
		Note that Theorem \ref{thm: DG Holder contin} shows that the function 
		\begin{align*}
			x\longmapsto\lm_x
		\end{align*}
		is $\bt_\al$--H\"older on $X$. Thus to see the claim we consider the following calculation.
		\begin{align*}
		&\absval{\tr_{0,\ta^{-n}(x)}^n(\ind_{\ta^{-n}(x)})(w)-\tr_{0,\ta^{-n}(x')}^n(\ind_{\ta^{-n}(x')})(w')}\\
		&\leq \absval{(\lm_{\ta^{-n}(x)}^n)^{-1} \sum_{z\in T_{x}^{-n}(w)}e^{S_n\phi_{\ta^{-n}(x)}(z)}-(\lm_{\ta^{-n}(x')}^n)^{-1} \sum_{z\in T_{x}^{-n}(w)}e^{S_n\phi_{\ta^{-n}(x)}(z)}} \\
		&\qquad+ \absval{(\lm_{\ta^{-n}(x')}^n)^{-1} \sum_{z\in T_{x}^{-n}(w)}e^{S_n\phi_{\ta^{-n}(x)}(z)}-(\lm_{\ta^{-n}(x')}^n)^{-1} \sum_{z'\in T_{x'}^{-n}(w')}e^{S_n\phi_{\ta^{-n}(x')}(z')}} \\
		&\leq \sum_{z\in T_{x}^{-n}(w)}\absval{e^{S_n\phi_{\ta^{-n}(x)}(z)}}\absval{(\lm_{\ta^{-n}(x)}^n)^{-1}-(\lm_{\ta^{-n}(x')}^n)^{-1}}\\
		&\qquad+\absval{\lm_{\ta^{-n}(x')}^n}^{-1}\absval{ \sum_{z\in T_{x}^{-n}(w)}e^{S_n\phi_{\ta^{-n}(x)}(z)}- \sum_{z'\in T_{x'}^{-n}(w')}e^{S_n\phi_{\ta^{-n}(x')}(z')}}\\
		&= \absval{\lm_{\ta^{-n}(x)}^n}^{-1}\sum_{z\in T_{x}^{-n}(w)}\absval{e^{S_n\phi_{\ta^{-n}(x)}(z)}}\absval{1-\frac{\lm_{\ta^{-n}(x)}^n}{\lm_{\ta^{-n}(x')}^n}}\\
		&\qquad+\absval{\lm_{\ta^{-n}(x')}^n}^{-1}\absval{\tr_{\ta^{-n}(x)}^n\ind_{\ta^{-n}(x)}(w)-\tr_{\ta^{-n}(x')}^n\ind_{\ta^{-n}(x')}(w')}\\\
		&= \absval{\lm_{\ta^{-n}(x')}^n}^{-1}\tr_{0,\ta^{-n}(x)}^n\ind_{\ta^{-n}(x)}(w)\absval{\lm_{\ta^{-n}(x)}^n-\lm_{\ta^{-n}(x')}^n}\\
		&\qquad+\absval{\lm_{\ta^{-n}(x')}^n}^{-1}\absval{\tr_{\ta^{-n}(x)}^n\ind_{\ta^{-n}(x)}(w)-\tr_{\ta^{-n}(x')}^n\ind_{\ta^{-n}(x')}(w')}\\
		&\leq C\cdot\absval{\lm_{\ta^{-n}(x')}^n}^{-1}\absval{\int_{\cJ_x}\tr_{\ta^{-n}(x)}^n\ind_{\ta^{-n}(x)}\,d\nu_x-\int_{\cJ_{x'}}\tr_{\ta^{-n}(x')}^n\ind_{\ta^{-n}(x')}\,d\nu_{x'}}\\
		&\qquad+C\cdot\absval{\lm_{\ta^{-n}(x')}^n}^{-1}\cdot\tr_{\ta^{-n}(x')}^n\ind_{\ta^{-n}(x')}(w')\cdot d_{\cJ}^\al((x,w),(x',w'))\\
		&\leq C\cdot\absval{\lm_{\ta^{-n}(x')}^n}^{-1}\tr_{\ta^{-n}(x')}^n(\ind_{\ta^{-n}(x')})(w)\cdot d^{\bt_\al}(x,x')\\
		&\qquad+C\cdot\tr_{0,\ta^{-n}(x')}^n\ind_{\ta^{-n}(x')}(w')\cdot d_{\cJ}^\al((x,w),(x',w'))\\
		&=C\cdot\tr_{0,\ta^{-n}(x')}^n(\ind_{\ta^{-n}(x')})(w)\cdot d^{\bt_\al}(x,x')\\
		&\qquad+C\cdot\tr_{0,\ta^{-n}(x')}^n\ind_{\ta^{-n}(x')}(w')\cdot d_{\cJ}^\al((x,w),(x',w'))\\
		&\leq C\cdot d_{\cJ}^{\vkp}((x,w),(x',w')),
		\end{align*}
		where $\vkp:=\min\set{\al,\bt_\al}$. The proof is now complete. 
\end{proof}

We now wish to show the function 	
\begin{align*}
	x\longmapsto\int_{\cJ_x}\tr_{r_{n-1},\ta^{-1}(x)}\circ\cdots\circ\tr_{r_0,\ta^{-n}(x)}(u_{\ta^{-n}(x)})\,d\nu_x
\end{align*}  
is H\"older continuous on $X$ by showing the integrand is H\"older continuous on $\cJ$ and then applying the Denker--Gordin H\"older continuity theorem (Theorem \ref{thm: DG Holder contin}).
\begin{align*}
&\absval{\tr_{r_{n-1},\ta^{-1}(x)}\circ\cdots\circ\tr_{r_0,\ta^{-n}(x)}(u_{\ta^{-n}(x)})(w) -\tr_{r_{n-1},\ta^{-1}(x')}\circ\cdots\circ\tr_{r_0,\ta^{-n}(x')}(u_{\ta^{-n}(x')})(w') }\\
&=\absval{\tr_{0,\ta^{-n}(x)}^n(e^{i\ol{S}_ng_{\ta^{-n}(x)}}u_{\ta^{-n}(x)})(w) -\tr_{0,\ta^{-n}(x')}^n(e^{i\ol{S}_ng_{\ta^{-n}(x')}}u_{\ta^{-n}(x')})(w') }\\
&=\left|(\lm_{\ta^{-n}(x)}^n)^{-1}\sum_{z\in T_{x}^{-n}(w)}e^{S_n\phi_{\ta^{-n}(x)}(z)}e^{i\ol{S}_ng_{\ta^{-n}(x)}(z)}u_{\ta^{-n}(x)}(z) \right.\\
&\qquad\qquad \left. - (\lm_{\ta^{-n}(x')}^n)^{-1}\sum_{z'\in T_{x'}^{-n}(w')}e^{S_n\phi_{\ta^{-n}(x')}(z')}e^{i\ol{S}_ng_{\ta^{-n}(x')}(z')}u_{\ta^{-n}(x')}(z')\right|.
\end{align*}
In light of \eqref{frm:tran index} we may rewrite the last equality from above as 
\begin{align*}
&\absval{\tr_{r_{n-1},\ta^{-1}(x)}\circ\cdots\circ\tr_{r_0,\ta^{-n}(x)}(u_{\ta^{-n}(x)})(w) -\tr_{r_{n-1},\ta^{-1}(x')}\circ\cdots\circ\tr_{r_0,\ta^{-n}(x')}(u_{\ta^{-n}(x')})(w') }\\
&=\left|(\lm_{\ta^{-n}(x)}^n)^{-1}\sum_{z\in T_{x}^{-n}(w)}e^{S_n\phi_{\ta^{-n}(x)}(z)}e^{i\ol{S}_ng_{\ta^{-n}(x)}(z)}u_{\ta^{-n}(x)}(z) \right.\\
&\qquad\qquad \left. - (\lm_{\ta^{-n}(x')}^n)^{-1}\sum_{z\in T_{x}^{-n}(w)}e^{S_n\phi_{\ta^{-n}(x')}(z')}e^{i\ol{S}_ng_{\ta^{-n}(x')}(z')}u_{\ta^{-n}(x')}(z')\right|
\end{align*}
We can then split this difference in to the sum of four differences in the standard way, which we call $(\Dl_1),\dots,(\Dl_4)$, that is 
\begin{align*}
&\left|(\lm_{\ta^{-n}(x)}^n)^{-1}\sum_{z\in T_{x}^{-n}(w)}e^{S_n\phi_{\ta^{-n}(x)}(z)}e^{i\ol{S}_ng_{\ta^{-n}(x)}(z)}u_{\ta^{-n}(x)}(z) \right.\\
&\qquad\qquad \left. - (\lm_{\ta^{-n}(x')}^n)^{-1}\sum_{z\in T_{x}^{-n}(w)}e^{S_n\phi_{\ta^{-n}(x')}(z')}e^{i\ol{S}_ng_{\ta^{-n}(x')}(z')}u_{\ta^{-n}(x')}(z')\right|\\
&\leq\left|(\lm_{\ta^{-n}(x)}^n)^{-1}\sum_{z\in T_{x}^{-n}(w)}e^{S_n\phi_{\ta^{-n}(x)}(z)}e^{i\ol{S}_ng_{\ta^{-n}(x)}(z)}u_{\ta^{-n}(x)}(z) \right.\\
&\qquad\qquad \left. - (\lm_{\ta^{-n}(x')}^n)^{-1}\sum_{z\in T_{x}^{-n}(w)}e^{S_n\phi_{\ta^{-n}(x)}(z)}e^{i\ol{S}_ng_{\ta^{-n}(x)}(z)}u_{\ta^{-n}(x)}(z)\right|\qquad& (\Dl_1)
\\  
&+\left|(\lm_{\ta^{-n}(x')}^n)^{-1}\sum_{z\in T_{x}^{-n}(w)}e^{S_n\phi_{\ta^{-n}(x)}(z)}e^{i\ol{S}_ng_{\ta^{-n}(x)}(z)}u_{\ta^{-n}(x)}(z) \right.&\\
&\qquad\qquad \left. - (\lm_{\ta^{-n}(x')}^n)^{-1}\sum_{z\in T_{x}^{-n}(w)}e^{i\ol{S}_ng_{\ta^{-n}(x)}(z)}u_{\ta^{-n}(x)}(z)\cdot e^{S_n\phi_{\ta^{-n}(x')}(z')}\right|\qquad& (\Dl_2)\\  	
&+\left|(\lm_{\ta^{-n}(x')}^n)^{-1}\sum_{z\in T_{x}^{-n}(w)}e^{i\ol{S}_ng_{\ta^{-n}(x)}(z)}u_{\ta^{-n}(x)}(z)\cdot e^{S_n\phi_{\ta^{-n}(x')}(z')} \right.&\\
&\qquad\qquad \left. - (\lm_{\ta^{-n}(x')}^n)^{-1}\sum_{z\in T_{x}^{-n}(w)}u_{\ta^{-n}(x)}(z)\cdot e^{S_n\phi_{\ta^{-n}(x')}(z')}e^{i\ol{S}_ng_{\ta^{-n}(x')}(z')}\right|\qquad &(\Dl_3)\\ 
&+\left|(\lm_{\ta^{-n}(x')}^n)^{-1}\sum_{z\in T_{x}^{-n}(w)}u_{\ta^{-n}(x)}(z)\cdot e^{S_n\phi_{\ta^{-n}(x')}(z')}e^{i\ol{S}_ng_{\ta^{-n}(x')}(z')}\right.&\\
&\qquad\qquad \left. - (\lm_{\ta^{-n}(x')}^n)^{-1}\sum_{z\in T_{x}^{-n}(w)} e^{S_n\phi_{\ta^{-n}(x')}(z')}e^{i\ol{S}_ng_{\ta^{-n}(x')}(z')}u_{\ta^{-n}(x')}(z')\right|\qquad& (\Dl_4)\\ 
\end{align*}
We now estimate each of the previous differences $(\Dl_1)$--$(\Dl_4)$, beginning with $(\Dl_1)$. Theorem \ref{thm: DG Holder contin} and Lemma \ref{lem: tr op HC} allows us to write
\begin{align*}
(\Dl_1)&=\absval{\sum_{z\in T_x^{-n}(w)}e^{S_n\phi_{\ta^{-n}(x)}(z)}e^{i\ol{S}_ng_{\ta^{-n}(x)}(z)}u_{\ta^{-n}(x)}(z)}\cdot\absval{(\lm_{\ta^{-n}(x)}^n)^{-1}-(\lm_{\ta^{-n}(x')}^n)^{-1}} \\
&=\absval{\lm_{\ta^{-n}(x)}^n}\absval{\tr_{0,\ta^{-n}(x)}^n(e^{i\ol{S}_ng_{\ta^{-n}(x)}}u_{\ta^{-n}(x)})(w)}\absval{(\lm_{\ta^{-n}(x)}^n)^{-1}-(\lm_{\ta^{-n}(x')}^n)^{-1}} \\
&\leq\tr_{0,\ta^{-n}(x)}^n(\absval{u_{\ta^{-n}(x)}})(w)\cdot\absval{1-\frac{\lm_{\ta^{-n}(x)}^n}{\lm_{\ta^{-n}(x')}^n}	}\\
&\leq C\cdot\absval{u_x(w)}\cdot\absval{\lm_{\ta^{-n}(x)}^n}^{-1}\absval{\lm_{\ta^{-n}(x)}^n-\lm_{\ta^{-n}(x')}^n}\\
&\leq C\cdot\absval{\lm_{\ta^{-n}(x)}^n}^{-1}\absval{\int_{\cJ_x}\tr_{\ta^{-n}(x)}^n\ind_{\ta^{-n}(x)}\,d\nu_x-\int_{\cJ_{x'}}\tr_{\ta^{-n}(x')}^n\ind_{\ta^{-n}(x')}\,d\nu_{x'}}\\
&\leq C\cdot\absval{\lm_{\ta^{-n}(x)}^n}^{-1}\cdot\tr_{\ta^{-n}(x')}^n(\ind_{\ta^{-n}(x')})(w)\cdot d^{\bt_\al}(x,x')\\
&=C\cdot\tr_{0,\ta^{-n}(x')}^n(\ind_{\ta^{-n}(x')})(w)\cdot d^{\bt_\al}(x,x')\\
&\leq C\cdot d^{\bt_\al}(x,x').
\end{align*}
Using \eqref{ineq: Sn DG*}, the difference $(\Dl_2)$ can be estimated as
\begin{align*}
(\Dl_2)&=\absval{(\lm_{\ta^{-n}(x')}^n)^{-1}\sum_{z\in T_{x}^{-n}(w)}e^{i\ol{S}_ng_{\ta^{-n}(x)}(z)}u_{\ta^{-n}(x)}(z)
	\left(e^{S_n\phi_{\ta^{-n}(x)}(z)}-e^{S_n\phi_{\ta^{-n}(x')}(z')}\right)}\\
&\leq \absval{\lm_{\ta^{-n}(x')}^n}^{-1}\sum_{z\in T_{x}^{-n}(w)}\absval{u_{\ta^{-n}(x)}(z)}\absval{e^{S_n\phi_{\ta^{-n}(x)}(z)}-e^{S_n\phi_{\ta^{-n}(x')}(z')}}\\
&\leq C\cdot	\absval{\lm_{\ta^{-n}(x')}^n}^{-1}\sum_{z\in T_{x}^{-n}(w)}e^{S_n\phi_{\ta^{-n}(x')}(z')}\absval{e^{S_n\phi_{\ta^{-n}(x)}(z)-S_n\phi_{\ta^{-n}(x')}(z')}-1}\\
&\leq C\cdot	\absval{\lm_{\ta^{-n}(x')}^n}^{-1}\sum_{z\in T_{x}^{-n}(w)}e^{S_n\phi_{\ta^{-n}(x')}(z')}e^{\absval{S_n\phi_{\ta^{-n}(x)}(z)-S_n\phi_{\ta^{-n}(x')}(z')}}\absval{S_n\phi_{\ta^{-n}(x)}(z)-S_n\phi_{\ta^{-n}(x')}(z')}\\
&\leq C\cdot d_{\cJ}^\al((x,w),(x',w'))\cdot	\absval{\lm_{\ta^{-n}(x')}^n}^{-1}\sum_{z\in T_{x}^{-n}(w)}e^{S_n\phi_{\ta^{-n}(x')}(z')}\\
&=C\cdot d_{\cJ}^\al((x,w),(x',w'))\cdot\tr_{0,\ta^{-n}(x')}^n(\ind_{\ta^{-n}(x')})(w')\\
&\leq C\cdot d_{\cJ}^\al((x,w),(x',w')).
\end{align*} 
Similarly, using \eqref{ineq: Snbar DG*}, the difference $(\Dl_3)$ can be estimated as 
\begin{align*}
(\Dl_3)&=\absval{(\lm_{\ta^{-n}(x')}^n)^{-1}\sum_{z\in T_{x}^{-n}(w)}u_{\ta^{-n}(x)}(z)
	\cdot e^{S_n\phi_{\ta^{-n}(x')}(z')}\left(e^{i\ol{S}_ng_{\ta^{-n}(x)}(z)}-e^{i\ol{S}_ng_{\ta^{-n}(x')}(z')}\right)}\\
&\leq\absval{\lm_{\ta^{-n}(x')}^n}^{-1}\sum_{z\in T_{x}^{-n}(w)}u_{\ta^{-n}(x)}(z)
\cdot e^{S_n\phi_{\ta^{-n}(x')}(z')}\absval{e^{i\ol{S}_ng_{\ta^{-n}(x)}(z)}-e^{i\ol{S}_ng_{\ta^{-n}(x')}(z')}}\\
&\leq C\cdot\absval{\lm_{\ta^{-n}(x')}^n}^{-1}\sum_{z\in T_{x}^{-n}(w)}
e^{S_n\phi_{\ta^{-n}(x')}(z')}\cdot\absval{\ol{S}_ng_{\ta^{-n}(x)}(z)-\ol{S}_ng_{\ta^{-n}(x')}(z')}\\
&\leq C\cdot d_{\cJ}^\al((x,w),(x',w'))\absval{\lm_{\ta^{-n}(x')}^n}^{-1}
\sum_{z\in T_{x}^{-n}(w)}e^{S_n\phi_{\ta^{-n}(x')}(z')}\\
&\leq C\cdot d_{\cJ}^\al((x,w),(x',w'))\cdot\tr_{0,\ta^{-n}(x')}^n(\ind_{\ta^{-n}(x')})(w')\\
&\leq C\cdot d_{\cJ}^\al((x,w),(x',w')).
\end{align*}
Given that $u\in\sH_\tau(\cJ)$, the final difference $(\Dl_4)$ can be estimated as 
\begin{align*}
(\Dl_4)&=\absval{(\lm_{\ta^{-n}(x')}^n)^{-1}\sum_{z\in T_{x}^{-n}(w)}e^{S_n\phi_{\ta^{-n}(x')}(z')}e^{i\ol{S}_ng_{\ta^{-n}(x')}(z')}\cdot \left(u_{\ta^{-n}(x)}(z)-u_{\ta^{-n}(x')}(z') \right) }\\
&\leq \absval{\lm_{\ta^{-n}(x')}^n}^{-1}\sum_{z\in T_{x}^{-n}(w)}e^{S_n\phi_{\ta^{-n}(x')}(z')}\cdot\absval{u_{\ta^{-n}(x)}(z)-u_{\ta^{-n}(x')}(z')}\\
&\leq C\cdot d_{\cJ}^\tau((x,w),(x',w'))\cdot\absval{\lm_{\ta^{-n}(x')}^n}^{-1}\sum_{z\in T_{x}^{-n}(w)}e^{S_n\phi_{\ta^{-n}(x')}(z')}\\
&=C\cdot d_{\cJ}^\tau((x,w),(x',w'))\cdot\tr_{0,\ta^{-n}(x')}^n(\ind_{\ta^{-n}(x')})(w')\\
&\leq C\cdot d_{\cJ}^\tau((x,w),(x',w')).
\end{align*}
All together this gives that 
\begin{align*}
	&\absval{\tr_{r_{n-1},\ta^{-1}(x)}\circ\cdots\circ\tr_{r_0,\ta^{-n}(x)}(u_{\ta^{-n}(x)})(w) -\tr_{r_{n-1},\ta^{-1}(x')}\circ\cdots\circ\tr_{r_0,\ta^{-n}(x')}(u_{\ta^{-n}(x')})(w') }\\
	&\leq \Dl_1+\Dl_2+\Dl_3+\Dl_4\\
	&\leq C\cdot d_{\cJ}^\al((x,w),(x',w'))+C\cdot d^{\bt_\al}(x,x')+C\cdot d_{\cJ}^\tau((x,w),(x',w'))\\
	&\leq C\cdot d_{\cJ}^\zt((x,w),(x',w')),	
\end{align*}
where 
\begin{align*}
	\zt=\min\set{\al,\bt_\al,\tau}.
\end{align*}
Upon application of Theorem \ref{thm: DG Holder contin} we see that the function 
\begin{align*}
x\longmapsto\int_{\cJ_x}\tr_{r_{n-1},\ta^{-1}(x)}\circ\cdots\circ\tr_{r_0,\ta^{-n}(x)}(u_{\ta^{-n}(x)})\,d\nu_x
\end{align*} 
is $\bt_\zt$--H\"older continuous on $X$ with uniformly bounded $\HH_{\bt_\zt}$ norm.
The proof of Proposition \ref{lem: perturbed int holder} is now complete.

\end{proof}
\begin{remark}
	The proof of Proposition \ref{lem: perturbed int holder} gives more. We have actually shown that the functions $\tr_{x}^n(u_x)$, $\tr_{0,x}^n(u_x)$, and $\tr_{r,x}^n(u_x)$ are each H\"older continuous on $\cJ$, which applying Theorem \ref{thm: DG Holder contin}, would mean that each of the functions  $x\longmapsto\tr_{x}^n(u_x)$, $x\longmapsto\tr_{0,x}^n(u_x)$, and $x\longmapsto\tr_{r,x}^n(u_x)$ are each H\"older continuous on $X$.
\end{remark}

\begin{remark}
	Considering Remark \ref{rem:DG* system ineq}, we see that Proposition \ref{lem: perturbed int holder}, provided that we know that the function $\rho$ is H\"older on $\cJ$, implies that the function
	\begin{align*}
	x\longmapsto\int_{\cJ_x}\tr_{r_{n-1},\ta^{-1}(x)}\circ\cdots\circ\tr_{r_0,\ta^{-n}(x)}(\rho_{\ta^{-n}(x)})\,d\nu_x
	\end{align*}  
	is H\"older continuous on $X$ with uniformly bounded $\HH$ norm, which satisfies condition \ref{assum 6} of Theorem \ref{thm: modified gouezel asip}. 
\end{remark}
In order to show the function $\rho$ is indeed H\"older continuous over $\cJ$ we will first need the following result while follows from the proof of Lemma 3.8 of \cite{mayer_distance_2011}.
\begin{lemma}
	For each $x\in X$ the sequence 
	\begin{align*}
		\rho_{x,n}:=\frac{1}{n}\sum_{k=0}^{n-1}\tr_{0,\ta^{-k}(x)}^k\ind_{\ta^{-k}(x)}
	\end{align*}
	is equicontinuous, i.e. there exists a subsequence $n_j$ such that 
	\begin{align*}
		\rho_{x,n_j}\lra\rho_x
	\end{align*} 
	uniformly.
\end{lemma}

Thus, in light of the previous remark, the following lemma establishes that the hypotheses of Theorem \ref{thm: modified gouezel asip} hold for DG*--systems.

\begin{lemma}
	The function $\rho$ is $\vkp$--H\"older on $X$, where $\vkp=\min\set{\al,\bt_\al}$.
\end{lemma}
\begin{proof}
	
	 To see this we appeal to the H\"older continuity of the function $\tr_{0,\ta^{-n}(x)}^n\ind_{\ta^{-n}(x)}$ which we just showed in Lemma \ref{lem: norm tr op HC}. Since we have that the functions 
	\begin{align*}
	\rho_{x,n}:=\frac{1}{n}\sum_{k=0}^{n-1}\tr_{0,\ta^{-k}(x)}^k\ind_{\ta^{-k}(x)}
	\end{align*}
	converge to $\rho_x$ for each $x\in X$ and from the proof of Lemma \ref{lem: norm tr op HC}, we also have that for $(x,w),(x',w')\in\cJ$ and each $n\in\NN$
	\begin{align*}
	\absval{\rho_{x,n}(w)-\rho_{x',n}(w')}&=\absval{\frac{1}{n}\sum_{k=0}^{n-1}\tr_{0,\ta^{-k}(x)}^k\ind_{\ta^{-k}(x)}-\frac{1}{n}\sum_{k=0}^{n-1}\tr_{0,\ta^{-k}(x')}^k\ind_{\ta^{-k}(x')}}\leq C\cdot d_{\cJ}^\vkp((x,w)(x',w')).
	\end{align*}
	Thus we have that 
	\begin{align*}
	\absval{\rho_{x}(w)-\rho_{x'}(w')}\leq C\cdot d_{\cJ}^\vkp((x,w),(x',w')).
	\end{align*}
	In particular, we see that for $z\in T_{x}^{-n}(w)$ and $z'\in T_{x'}^{-n}(w')$
	\begin{align*}
	\absval{\rho_{\ta^{-n}(x)}(z)-\rho_{\ta^{-n}(x)'}(z')}\leq C\cdot \gm^{-\vkp n}\cdot d_{\cJ}^\vkp((x,w),(x',w')).
	\end{align*}
	This finishes the proof.
\end{proof}
We have finally shown that the hypotheses of Theorem \ref{thm: modified gouezel asip} hold for DG*--systems, and we now have the following.
\begin{theorem}
	Let $T:\jl\lra\jl$ be a DG*--system and $g\in\sH_\al^*(\cJ)$. Then either there exists a number $\sg^2>0$ such that the process $\set{g\circ T^n-\mu(g)}_{n\in\NN}$ satisfies an ASIP with limiting covariance $\sg^2$ for any error exponent larger that $1/4$, or, if $\sg^2=0$, then we have that
	\begin{align*}
	\sup_{n\in\NN}\norm{S_ng-\mu(g)}_{L^2(\mu)}<\infty .
	\end{align*}
\end{theorem}

\section*{Acknowledgments}
The author would like to thank Mariusz Urba\'nski for the many thoughtful discussions which inspired this work.

%\newpage
%\nocite{*}
\bibliographystyle{jabbrv_abbrv}
\bibliography{mybib} 

\end{document}